\documentclass[10pt]{amsart}
\usepackage{amsfonts}
\usepackage{graphicx}
\usepackage{amscd}
\usepackage{amsmath,amsfonts,amssymb,amsthm,mathrsfs,xypic}
\xyoption{curve}
\usepackage{fancybox}
\newtheorem{thm}{Theorem}[section]

\newtheorem{cor}[thm]{Corollary}

\newtheorem{lem}[thm]{Lemma}

\newtheorem{prop}[thm]{Proposition}
\newtheorem{rem}[thm]{Remark}

\numberwithin{equation}{section}

\newcommand{\ik}{\operatorname {IKer}}

\newcommand{\s}{\hfill\blacksquare}

\newcommand{\End}{\operatorname{End}}
\newcommand{\add}{\operatorname{add}}
\newcommand{\Ima}{\operatorname{Im}}
\newcommand{\Cok}{\operatorname{Coker}}
\newcommand{\Ker}{\operatorname{Ker}}
\newcommand{\Hom}{\operatorname{Hom}}
\newcommand{\Ext}{\operatorname{Ext}}

\newcommand{\gld}{\operatorname{gl.dim}}
\begin{document}
\title [Monomorphism categories and cotilting theory]{Monomorphism categories, cotilting
theory,  and Gorenstein-projective modules}
\author [Pu Zhang] {Pu Zhang}
\thanks{Supported by the NSF of China (10725104), and STCSM (09XD1402500).}
\thanks{pzhang$\symbol{64}$sjtu.edu.cn}
\dedicatory {Dedicated to Claus Michael Ringel on the occasion of
his 65$^{th}$ birthday}

\maketitle

\centerline{Department of Mathematics, \ \ Shanghai Jiao Tong
University} \centerline{Shanghai 200240, P. R. China}
\begin{abstract}
\ The monomorphism category $\mathcal S_n(\mathcal X)$ is
introduced, where $\mathcal X$ is a full subcategory of the module
category $A$-mod of Artin algebra $A$. The key result is a
reciprocity of the monomorphism operator $\mathcal S_n$ and the left
perpendicular operator $^\perp$: for a cotilting $A$-module $T$,
there is a canonical construction of a cotilting $T_n(A)$-module
${\rm \bf m}(T)$, such that $\mathcal S_n(^\perp T) = \ ^\perp {\rm
\bf m}(T)$.

As applications, $\mathcal S_n(\mathcal X)$ is a resolving
contravariantly finite subcategory in $T_n(A)$-mod with
$\widehat{\mathcal S_n(\mathcal X)} = T_n(A)$-mod if and only if
$\mathcal X$ is a resolving contravariantly finite subcategory in
$A$-mod with $\widehat{\mathcal X} = A$-mod. For a Gorenstein
algebra $A$, the category $T_n(A)\mbox{-}\mathcal Gproj$ of
Gorenstein-projective $T_n(A)$-modules can be explicitly determined
as $\mathcal S_n(^\perp A)$. Also, self-injective algebras $A$ can
be characterized by the property $T_n(A)\mbox{-}\mathcal Gproj =
\mathcal S_n(A)$. Using  $\mathcal S_n(A)= \ ^\perp {\rm \bf
m}(D(A_A))$, a characterization of $\mathcal S_n(A)$ of finite type
is obtained.

\vskip5pt

\end{abstract}

\vskip10pt

\centerline {\bf  Introduction}

\vskip10pt

Throughout $A$ is an Artin algebra, and $n\ge 2$ an integer. Let
$A$-mod be the category of finitely generated left $A$-modules, and
$\mathcal X$ a full subcategory of $A$-mod. Denote by ${\rm
Mor}_n(A)$ the morphism category, which is equivalent to
$T_n(A)$-mod, where $T_n(A)$ is the upper triangular matrix algebra
of $A$. Let $\mathcal S_n(A)$ denote the full subcategory of ${\rm
Mor}_n(A)$ given by $\mathcal S_n(A) = \{X_{(\phi_i)}\in {\rm
Mor}_n(A) \ | \ \mbox{all} \ \phi_i
 \ \mbox{are monomorphisms} \}.$

\vskip5pt G. Birkhoff [Bir] initiated to classify the indecomposable
objects of $\mathcal S_2(\Bbb Z/\langle p^t\rangle)$. In [RW] the
indecomposable objects of $\mathcal S_2(\Bbb Z/\langle p^t\rangle)$
with $t\le 5$ were determined. In [Ar] $\mathcal S_n(R)$ was denoted
by $\mathcal C(n, R)$, where $R$ is a commutative uniserial artinian
ring; and the complete lists of $\mathcal C(n, R)$ of finite type,
and of the representation types of $\mathcal C(n, k[x]/\langle
x^t\rangle)$, have been given by D. Simson [S] (see also [SW]). C.
M. Ringel and M. Schmidmeier ([RS1] - [RS3]) have intensively
studied the monomorphism category $\mathcal S_2(A)$. According to
[RS2], $\mathcal S_n(A)$ is a functorially finite subcategory in
$T_n(A)$-mod; and hence $\mathcal S_n(A)$ has Auslander-Reiten
sequences. For more recent work related to the monomorphism
categories we refer to [C], [IKM] and [KLM].

\vskip5pt

On the other hand, M. Auslander and I. Reiten [AR] have established
a relation between resolving contravariantly finite subcategories
and cotilting theory, by asserting that $\mathcal X$ is resolving
and contravariantly finite with $\widehat{\mathcal X} = A$-mod if
and only if $\mathcal X = \ ^\perp T$ for some cotilting $A$-module
$T$ ([AR], Theorem 5.5(a)), where $^\perp T$ is the left
perpendicular category of $T$.

\vskip5pt

Define $\mathcal S_n(\mathcal X)$ to be the full subcategory of
${\rm Mor}_n(A)$ of all the objects $X_{(\phi_i)},$ where all
$X_i\in\mathcal X$, all $\phi_i$ are monomorphisms, and all
$\Cok\phi_i \in\mathcal X$. A main problem we concern is: when is
$\mathcal S_n(\mathcal X)$ contravariantly finite in $T_n(A)$-mod?
This leads to the following reciprocity of the monomorphism operator
$\mathcal S_n$ and the left perpendicular operator \ $^\perp$: given
a cotilting $A$-module $T$, then ${\rm \bf m}(T) =
\left(\begin{smallmatrix}
   T \\
   0\\ 0\\ \vdots \\ 0
\end{smallmatrix}\right) \oplus \left(\begin{smallmatrix}
  T \\
  T\\ 0 \\ \vdots \\ 0
\end{smallmatrix}\right)\oplus \cdots \oplus \left(\begin{smallmatrix}
  T  \\ T \\ T \\ \vdots \\ T
\end{smallmatrix}\right)$ is a cotilting $T_n(A)$-module, such
that $\mathcal S_n(^\perp T) = \ ^\perp {\rm \bf m}(T)$. See Theorem
3.1. The proof needs the contravariantly finiteness of $\mathcal
S_n(A)$ in $T_n(A)$-mod, and the six adjoint pairs between $A$-mod
and $T_n(A)$-mod.

\vskip10pt

We illustrate this reciprocity with several applications. First, we
have a solution to the main problem: $\mathcal S_n(\mathcal X)$ is a
resolving contravariantly finite subcategory in $T_n(A)$-mod with
$\widehat{\mathcal S_n(\mathcal X)} = T_n(A)$-mod if and only if
$\mathcal X$ is a resolving contravariantly finite subcategory in
$A$-mod with $\widehat{\mathcal X} = A$-mod (Theorem 3.9).

\vskip5pt

As another application, taking $T = \ _AA$ for Gorenstein algebra
$A$, the category $T_n(A)\mbox{-}\mathcal Gproj$ of
Gorenstein-projective $T_n(A)$-modules, can be  determined
explicitly as $\mathcal S_n(^\perp A)$ (Corollary 4.1). By D.
Happel's triangle-equivalence $D^b(A)/K^b(\mathcal P(A)) \cong
\underline {^\perp A}$ for Gorenstein algebra $A$, one has
$D^b(T_n(A))/K^b(\mathcal P(T_n(A))) = \underline {\mathcal
S_n(^\perp A)}$ (Corollary 4.3). Also, self-injective algebras $A$
can be characterized by the property $\mathcal S_n(A) =
T_n(A)\mbox{-}\mathcal Gproj$ (Theorem 4.4).

\vskip5pt

The representation type of $\mathcal S_n(A)$ is quite different from
the ones of $A$ and of $T_n(A)$. For example, $k[x]/\langle
x^t\rangle$ is of finite type, $T_2(k[x]/\langle x^t\rangle)$ is of
finite type if and only if $t\le 3$, but $\mathcal S_2(k[x]/\langle
x^t\rangle)$ is of finite type if and only if $t\le 5$, where $k$ is
an algebraically closed field. If $t > 6$ then $\mathcal
S_2(k[x]/\langle x^t\rangle)$ is of ``wild" type, while $\mathcal
S_2(k[x]/\langle x^6\rangle)$ is of ``tame" type ([S], Theorems 5.2
and 5.5). A complete classification of indecomposable objects of
$\mathcal S_2(k[x]/\langle x^6\rangle)$ is exhibited in [RS3].
Inspired by Auslander's classical result: $A$ is of finite type if
and only if there is an $A$-generator-cogenerator $M$ such that
$\gld \End_A(M)$ $\le 2$ ([Au], Chapter III), by using  $\mathcal
S_n(A)= \ ^\perp {\rm \bf m}(D(A_A))$, we prove that $\mathcal
S_n(A)$ is of finite type if and only if there is a bi-generator $M$
of $\mathcal S_n(A)$ such that $\gld \End_{T_n(A)}(M)\le 2$ (Theorem
5.1). As a corollary, for a self-injective algebra $A$, $T_n(A)$ is
CM-finite if and only if there is a $T_n(A)$-generator $M$ which is
Gorenstein-projective, such that $\gld {\rm End}_{T_n(A)}(M) \leq 2$
(Corollary 5.2).
\section{\bf \ Monomorphism categories}
We will define the monomorphism category $\mathcal S_n(\mathcal X)$
and give its basic properties needed later.
\subsection{} \  An object of the morphism category
${\rm Mor}_n(A)$ is $X_{(\phi_i)}=\left(\begin{smallmatrix}
X_1\\
\vdots\\ X_n
\end{smallmatrix}\right)_{(\phi_i)}$, where $\phi_i: X_{i+1}\rightarrow X_i$
\vskip5pt \noindent is an $A$-map, \ $1\le i\le n-1$; and a morphism
$X_{(\phi_i)}\rightarrow Y_{(\theta_i)}$ is
$f=\left(\begin{smallmatrix} f_1\\
\vdots\\ f_n
\end{smallmatrix}\right)$,
where $f_i: X_i \rightarrow Y_i$ is an $A$-map, \ $1\le i\le n$,
such that every square in the following diagram commutes
\[\xymatrix {
X_n\ar[d]_{f_n}\ar[r]^{\phi_{n-1}} & X_{n-1}\ar[d]_{f_{n-1}}\ar[r]^{\phi_{n-2}} & \cdots \ar[r]^{\phi_{1}} & X_1\ar[d]_{f_1}\\
Y_n\ar[r]^-{\theta_{n-1}} & Y_{n-1}\ar[r]^{\theta_{n-2}} & \cdots
\ar[r]^{\theta_{1}} & Y_1.}\eqno(1.1)\] Note that ${\rm Mor}_n(A)$
is again an abelian category, and that a sequence
$Z_{(\psi_i)}\stackrel f\longrightarrow Y_{(\theta_i)}\stackrel g
\longrightarrow X_{(\phi_i)}$ in ${\rm Mor}_n(A)$ is exact  if and
only if $Z_i\stackrel {f_i}\longrightarrow Y_i\stackrel {g_i}
\longrightarrow X_i$ is exact in $A$-mod for each $1\le i\le n$.
\subsection \ We define $\mathcal S_n(\mathcal X)$ to be the full subcategory of ${\rm
Mor}_n(A)$ consisting of all the objects $X_{(\phi_i)},$ where
$X_i\in\mathcal X$ for $1\le i\le n$, \ $\phi_i:
X_{i+1}\hookrightarrow X_{i}$ is a monomorphism and
$\Cok\phi_i\in\mathcal X$ for $1\le i\le n-1$. In particular, we
have  $\mathcal S_n(A\mbox{-}{\rm mod}) = \{X_{(\phi_i)}\in {\rm
Mor}_n(A) \ | \ \phi_i
 \ \mbox{is monic,} \ 1\le i\le n-1\},$ which will be denoted by $\mathcal
 S_n(A)$. Dually, $\mathcal F_n(\mathcal X)$ is the full subcategory of ${\rm
Mor}_n(A)$ consisting of all the objects $X_{(\phi_i)},$ where
$X_i\in\mathcal X, \ 1\le i\le n$, \ $\phi_i: \
X_{i+1}\twoheadrightarrow X_{i}$ is an epimorphism and
$\Ker\phi_i\in\mathcal X$ for $1\le i\le n-1$. We call $\mathcal
S_n(\mathcal X)$ and $\mathcal F_n(\mathcal X)$ {\it the
monomorphism category} and {\it the epimorphism category of}
$\mathcal X$, respectively.

\vskip10pt\begin{lem} \ Let $A$ be an Artin algebra and $\mathcal X$
a full subcategory of $A$-mod.

\vskip5pt

$(i)$ \ Let $0\rightarrow Z_{(\psi_i)} \rightarrow Y_{(\theta_i)}
\rightarrow X_{(\phi_i)} \rightarrow 0$ be an exact sequence in
${\rm Mor}_n(A)$. Then the following induced sequences are exact for
each $1\le i\le n-1$
\begin{align*}0&\longrightarrow \Ker (\psi_{1}\cdots \psi_{i})
\longrightarrow \Ker (\theta_{1}\cdots \theta_{i})\longrightarrow
\Ker (\phi_{1}\cdots \phi_{i})\longrightarrow
\\ &\longrightarrow \Cok (\psi_{1}\cdots \psi_{i})\longrightarrow \Cok
(\theta_{1}\cdots \theta_{i})\longrightarrow \Cok (\phi_{1}\cdots
\phi_{i})\longrightarrow 0,\end{align*} and $$0\longrightarrow \Ker
\psi_{i} \longrightarrow \Ker \theta_{i}\longrightarrow \Ker
\phi_{i}\longrightarrow \Cok \psi_{i}\longrightarrow \Cok
\theta_{i}\longrightarrow \Cok \phi_{i}\longrightarrow 0.$$
\vskip5pt $(ii)$ \ $\mathcal S_n(\mathcal X)$  is closed under
extensions (resp., kernels of epimorphisms, direct summands) if and
only if $\mathcal X$ is closed under extensions (resp., kernels of
epimorphisms, direct summands).

\vskip5pt

$(iii)$ \ $\mathcal S_n(A)$ is closed under subobjects.

\vskip5pt

$(iv)$ \ If $\mathcal X$ is closed under extensions, then there is
an equivalence of categories $\mathcal S_n(\mathcal X)\cong \mathcal
F_n(\mathcal X)$ given by
$$\mathcal S_n(\mathcal X)\ni\left(\begin{smallmatrix}
X_1\\ X_2 \\ \vdots\\ X_{n-1} \\ X_n
\end{smallmatrix}\right)_{(\phi_i)} \mapsto \left(\begin{smallmatrix}
\Cok\phi_1\\ \Cok(\phi_1\phi_2) \\ \vdots\\ \Cok(\phi_1\cdots\phi_{n-1}) \\
X_1
\end{smallmatrix}\right)_{(\phi'_i)}$$
where $\phi'_i: \Cok (\phi_1\cdots \phi_{i+1})\twoheadrightarrow
\Cok (\phi_1\cdots \phi_{i}), \ 1\le i\le n-2,$ and $\phi'_{n-1}:
X_1\twoheadrightarrow \Cok(\phi_1\cdots \phi_{n-1})$, are the
canonical epimorphisms, with a quasi-inverse
$$\mathcal F_n(\mathcal X)\ni\left(\begin{smallmatrix}
X_1\\ X_2 \\ \vdots\\ X_{n-1} \\ X_n
\end{smallmatrix}\right)_{(\phi_i)} \mapsto \left(\begin{smallmatrix}
X_n\\ \Ker(\phi_1\cdots\phi_{n-1}) \\ \vdots\\ \Ker(\phi_{n-2}\phi_{n-1}) \\
\Ker\phi_{n-1}
\end{smallmatrix}\right)_{(\phi''_i)},$$
where $\phi''_i: \Ker(\phi_i\cdots\phi_{n-1})\hookrightarrow
\Ker(\phi_{i-1}\cdots\phi_{n-1}), \ 2\le i\le n-1,$ and $\phi''_1:
\Ker(\phi_1\cdots\phi_{n-1})\hookrightarrow X_n$, are the canonical
monomorphisms.
\end{lem}
\noindent{\bf Proof.} \ Applying Snake Lemma to the following
commutative diagram with exact rows
\[\xymatrix{0\ar[r] & Z_{i+1} \ar[r]\ar[d]_-{\psi_{1}\cdots \psi_{i}}
& Y_{i+1}\ar[r]\ar[d]^-{\theta_{1}\cdots \theta_{i}} &
X_{i+1}\ar[r]\ar[d]^-{\phi_{1}\cdots \phi_{i}} & 0
\\ 0\ar[r] & Z_{1} \ar[r] & Y_{1}\ar[r]& X_{1}\ar[r] & 0}\]
we get the first exact sequence in $(i)$; and the second one can be
similarly obtained. $(ii)$ follows from $(i)$; $(iii)$ can seen from
$(1.1)$, and $(iv)$ is clear. $\s$

\vskip10pt

\subsection{} \ Let $X = \left(\begin{smallmatrix}
X_1\\
\vdots\\ X_n
\end{smallmatrix}\right)_{(\phi_i)}\in {\rm Mor}_n(A).$
We call $X_i$ the $i$-th branch of $X$, and $\phi_i$ the $i$-th
morphism of $X$. For each $1\le i\le n$, we define a functor ${\rm
\bf m}_i: A\mbox{-}{\rm mod} \rightarrow \mathcal S_n(A)$ as
follows. For $M\in A$-mod, the $j$-th branch of ${\rm \bf m}_i(M)$
is $M$ if $j\le i$, and $0$ if $j>i$; and the $j$-th morphism of
${\rm \bf m}_i(M)$ is ${\rm id}_M$ if $j < i$, and $0$ if $j\ge i$.
For each $A$-map $f: M\rightarrow N$, we define
$${\rm \bf m}_i(f) = \left(\begin{smallmatrix}
f\\
\vdots\\ f \\ 0 \\ \vdots \\ 0
\end{smallmatrix}\right): \ \ {\rm \bf m}_i(M) = \left(\begin{smallmatrix}
M\\
\vdots\\ M \\ 0 \\ \vdots \\ 0
\end{smallmatrix}\right)\longrightarrow  {\rm \bf m}_i(N) = \left(\begin{smallmatrix}
N\\
\vdots\\ N \\ 0 \\ \vdots \\ 0
\end{smallmatrix}\right).$$
Note that the restriction of ${\rm \bf m}_i$ to $\mathcal X$ gives
a functor $\mathcal X \rightarrow \mathcal S_n(\mathcal X)$. A
functor ${\rm \bf p}_i: A\mbox{-}{\rm mod} \rightarrow \mathcal
F_n(A)$ is dually defined, $1\le i\le n$. The $j$-th branch of ${\rm
\bf p}_i(M)$ is $M$ if $j\ge n-i+1$, and $0$ if $j < n-i+1$; and the
$j$-th morphism of ${\rm \bf p}_i(M)$ is ${\rm id}_M$ if $j \ge n-
i+1$, and $0$ if $j < n- i+1$. Also we define
$${\rm \bf p}_i(f) = \left(\begin{smallmatrix}
0\\
\vdots\\ 0 \\ f \\ \vdots \\ f
\end{smallmatrix}\right): \ \ {\rm \bf p}_i(M) = \left(\begin{smallmatrix}
0\\
\vdots\\ 0 \\ M \\ \vdots \\ M
\end{smallmatrix}\right)\longrightarrow  {\rm \bf p}_i(N) = \left(\begin{smallmatrix}
0\\
\vdots\\ 0 \\ N \\ \vdots \\ N
\end{smallmatrix}\right).$$
The restriction of ${\rm \bf p}_i$ to $\mathcal X$ gives a functor
$\mathcal X \rightarrow \mathcal F_n(\mathcal X)$. We have ${\rm \bf
m}_n(M) = {\rm \bf p}_n(M), \ \forall \ M\in A$-mod.

\vskip10pt The following facts imply that in fact there are six
adjoint pairs between $A$-mod and ${\rm Mor}_n(A)$.

\vskip10pt
\begin{lem} \ Let $A$ be an Artin algebra.
Then for each object $X=X_{(\phi_i)}\in {\rm Mor}_n(A)$ and each
$A$-module $M$, we have isomorphisms of abelian groups, which are
natural in both positions
$$\Hom_{{\rm Mor}_n(A)}({\rm \bf m}_i(M), X) \cong \Hom_A (M, X_i), \  \ \ 1\le i\le n, \eqno(1.2)$$
$$\Hom_{{\rm Mor}_n(A)}(X, {\rm \bf m}_i(M)) \cong \Hom_A (\Cok(\phi_{1}\cdots \phi_{i}), M),
\ \  \ 1\le i\le n-1, \eqno(1.3)$$
$$\Hom_{{\rm
Mor}_n(A)}(X, {\rm \bf m}_n(M)) \cong \Hom_A(X_1, M), \eqno(1.4)$$
$$\Hom_{{\rm Mor}_n(A)}(X, {\rm \bf p}_i(M)) \cong \Hom_A (X_{n-i+1},
M), \  \ \ 1\le i\le n, \eqno(1.5)$$
$$\Hom_{{\rm Mor}_n(A)}({\rm \bf p}_i(M), X) \cong \Hom_A(M, \Ker(\phi_{n-i}\cdots
\phi_{n-1})), \ \ \ 1\le i\le n-1, \eqno(1.6)$$
$$\Hom_{{\rm Mor}_n(A)}({\rm \bf p}_n(M), X) \cong \Hom_A (M, X_n). \eqno(1.7)$$
\end{lem} \noindent{\bf Proof.} \ We justify $(1.3)$. Let
$\pi_{i}: X_1 \twoheadrightarrow \Cok(\phi_{1}\cdots \phi_{i})$ be
the canonical epimorphism, $1\le i\le n-1$. Consider the
homomorphism of abelian groups $\Hom_A (\Cok(\phi_{1}\cdots
\phi_{i}), M)\rightarrow \Hom_{{\rm Mor}_n(A)}(X, {\rm \bf m}_i(M))$
given
by $$g \mapsto \left(\begin{smallmatrix} g\pi_{i}\\ g\pi_i\phi_1 \\ \vdots \\
g\pi_{i}\phi_{1}\cdots\phi_{i-1} \\ 0
\\ \vdots \\ 0
\end{smallmatrix}\right): \ \ X\longrightarrow \left(\begin{smallmatrix}
M\\ M\\ \vdots\\ M \\ 0 \\ \vdots \\ 0
\end{smallmatrix}\right), \ \forall \ g\in \Hom_A (\Cok(\phi_{1}\cdots \phi_{i}), M).$$
By $(1.1)$ we infer that it is surjective, and it is injective since
$\pi$ is epic. It is clear that the isomorphisms are natural in both
positions. $\s$
\subsection{} \  Let $\mathcal P(A)$ (resp. $\mathcal I(A)$) be
the full subcategory of $A$-mod of projective (resp. injective)
$A$-modules,  and ${\rm Ind}\mathcal P(A)$ (resp. ${\rm Ind}\mathcal
I(A)$) be the set of pairwise non-isomorphic indecomposable
projective (resp. injective) $A$-modules.

\vskip10pt
\begin{lem} \ Let $A$ be an Artin algebra. Then

\vskip5pt

$(i)$ \ There is an equivalence of categories ${\rm Mor}_n(A)\cong
T_n(A)\mbox{-}{\rm mod},$ which preserves the exact structures,
where $T_n(A)$ is the $n\times n$ upper triangular matrix algebra
$\left(\begin{smallmatrix}
A&A&\cdots &A\\
0&A&\cdots&A\\
& &\ddots &\\
\\ 0&0&\cdots&A
\end{smallmatrix}\right).$

$(ii)$ \ Under this equivalence,  we have
$${\rm Ind}\mathcal P(T_n(A)) = \{{\rm \bf m}_1(P), \  \cdots,  \ {\rm \bf m}_n(P) \ | \ P \in {\rm Ind}\mathcal P(A)\}
\subseteq \mathcal S _n(A), \eqno(1.8)$$
$${\rm Ind}\mathcal I(T_n(A)) =
\{{\rm \bf p}_1(I), \ \cdots, \ {\rm \bf p}_n(I) \ | \ I \in {\rm
Ind}\mathcal I(A)\}\subseteq \mathcal F_n(A). \eqno(1.9)$$
\end{lem}
\noindent {\bf Proof.} \ $(i)$ \  This is well-known, at least for
$n=2$ (see [ARS], p.71). For convenience we include a short
justification. For $1\le i\le j\le n$, let $e_{ij}\in T_n(A)$ be the
matrix with $1$ in the $(i,j)$-entry, and $0$ elsewhere. For a
$T_n(A)$-module $M$ we have $M = e_{11}M\oplus\cdots \oplus e_{nn}M$
as $A$-modules, and for an $A$-map $f: M\rightarrow N$, the
restriction $f_i$ of $f$ to $e_{ii}M$ gives an $A$-map $f_i: e_{ii}M
\rightarrow e_{ii}N$. Consider a functor $F: T_n(A)\mbox{-}{\rm mod}
\rightarrow {\rm Mor}_n(A)$ defined by $F(M) =
\left(\begin{smallmatrix}
   e_{11}M \\ \vdots \\ e_{nn}M
\end{smallmatrix}\right)_{(\phi_{M, i})}$,
where $\phi_{M,i}: e_{i+1 i+1}M\rightarrow e_{ii}M$ is the $A$-map
given by $\phi_{M, i}(e_{i+1 i+1}x) = e_{i i}e_{i i+1} e_{i+1 i+1}x
\in e_{i i}M, \ 1\le i\le n-1,$ and
$$F(f) = \left(\begin{smallmatrix}
   f_{1} \\ \vdots \\ f_{n}
\end{smallmatrix}\right): \left(\begin{smallmatrix}
   e_{11}M \\ \vdots \\ e_{nn}M
\end{smallmatrix}\right)_{(\phi_{M, i})} \longrightarrow \left(\begin{smallmatrix}
   e_{11}N \\ \vdots  \\ e_{nn}N
\end{smallmatrix}\right)_{(\phi_{N, i})}.$$
Then $F$ is fully faithful. For each object
$\left(\begin{smallmatrix}
X_1\\
\vdots\\
X_n
\end{smallmatrix}\right)_{(\phi_i)}\in {\rm Mor}_n(A),$ put $X =
\bigoplus\limits_{1\le i\le n}X_i$, and write an element of $X$ as
$\left(\begin{smallmatrix}
x_1\\
\vdots\\
x_n
\end{smallmatrix}\right)$ with $x_i\in X_i$. With  a $T_n(A)$-action on $X$ defined by
$$\left(\begin{smallmatrix}
a_{11}&a_{12}&\cdots &a_{1n}\\
0&a_{22}&\cdots&a_{2n}\\
& &\ddots &\\
\\ 0&0&\cdots&a_{nn}
\end{smallmatrix}\right)\left(\begin{smallmatrix}
x_1\\ x_2\\
\vdots\\
x_n
\end{smallmatrix}\right) = \left(\begin{smallmatrix}
a_{11}x_1 + \sum\limits_{j>1}a_{1j} \phi_1\cdots\phi_{j-1}(x_j)\\
\vdots \\ a_{ii}x_i + \sum\limits_{j>i}a_{ij}
\phi_i\cdots\phi_{j-1}(x_j) \\
\vdots\\
a_{nn}x_n
\end{smallmatrix}\right),$$
$X$ is a $T_n(A)$-module such that $F(X) = \left(\begin{smallmatrix}
X_1\\
\vdots\\
X_n
\end{smallmatrix}\right)_{(\phi_i)}$, i.e., $F$ is dense.

\vskip5pt

$(ii)$ \ Since $\End_{T_n(A)}({\rm \bf m}_i(M)) \cong \End_A(M)
\cong \End_{T_n(A)}({\rm \bf p}_i(M)), \ \forall \ M\in A$-mod, it
follows that if $M$ is indecomposable then ${\rm \bf m}_i(M)$ and
${\rm \bf p}_i(M)$ are indecomposable. Since ${\rm \bf m}_i$ and
${\rm \bf p}_i$ are additive functors, $(1.8)$ follows from the
decomposition $T_n(A) = \bigoplus\limits_{1\le i\le n} {\rm \bf
m}_i(A)$ as left $T_n(A)$-modules. By $(1.5)$ we see that ${\rm \bf
p}_i(I), \ 1\le i\le n$, are indecomposable injective
$T_n(A)$-modules, where $I\in {\rm Ind}\mathcal I(A)$, and then we
infer $(1.9)$, by comparing the number of pairwise non-isomorphic
indecomposable injective $T_n(A)$-modules. $\s$

\vskip10pt

From now on we identify $T_n(A)\mbox{-}{\rm mod}$ with ${\rm
Mor}_n(A)$.

\vskip10pt

\subsection{} \ A full subcategory $\mathcal X$ of $A$-mod is {\it resolving}
if $\mathcal X$ contains all projective $A$-modules, $\mathcal X$ is
closed under extensions, kernels of epimorphisms, and direct
summands. By Lemmas 1.3$(ii)$ and 1.1$(ii)$, and using functor ${\rm
\bf m}_1: \mathcal X\rightarrow \mathcal S_n(\mathcal X)$ we get

\vskip10pt

\begin{cor} \ Let $A$ be an Artin algebra and $\mathcal X$ a full subcategory of $A$-mod.
Then $\mathcal S_n(\mathcal X)$ is a resolving subcategory of ${\rm
Mor}_n(A)$ if and only if $\mathcal X$ is a resolving subcategory of
$A$-mod. \end{cor}

\section {\bf Functorially finiteness of $\mathcal S_n(A)$ in ${\rm Mor}_n(A)$}

The idea of the following result  comes from [RS2] for $\mathcal
S_2(A)$.

\vskip10pt

\begin{thm} \ (Ringel - Schmidmeier) \ Let $A$ be an Artin algebra. Then $\mathcal S_n(A)$
is a functorially finite subcategory in ${\rm Mor}_n(A)$ and has
Auslander-Reiten sequences.
\end{thm}

\subsection{} \ Let $X_{(\phi_i)}\in {\rm Mor}_n(A)$. Fix an injective envelope
$e_i': \Ker\phi_i\hookrightarrow \ik\phi_i.$ Define object ${\rm
rMon}(X)\in \mathcal S_n(A)$ as follows. We have an $A$-map $e_i:
X_{i+1}\rightarrow \ik\phi_i$ such that the following diagram
commutes for each \ $1\le i\le n-1$
\[\xymatrix{\Ker\phi_i \ar @{^{(}->}[r]\ar @{^{(}->}[d]_{e'_i} &
X_{i+1} \ar @{>} [ld]^{e_i}\\
\ik\phi_i  & &.}\eqno(2.1)
\] Of course $e_1, \cdots, e_{n-1}$ are not unique. However we
choose and fix them, and then define

\noindent $\theta_i:  X_{i+1}\oplus \ik\phi_{i+1}\oplus
\cdots\oplus\ik\phi_{n-1} \longrightarrow X_{i}\oplus \ik
\phi_i\oplus \ik\phi_{i+1}\oplus\cdots\oplus\ik\phi_{n-1}$  to be
$$\theta_i = \left(\begin{smallmatrix}
\phi_i& 0 & 0  &\cdots &  0 \\
e_i& 0 & 0 & \cdots &  0 \\
0 & 1& 0 & \cdots &  0 \\
0 & 0& 1 & \cdots &  0\\
 \vdots & \vdots & \vdots &
\cdots & \vdots \\
0 & 0& 0 & \cdots & 1
\end{smallmatrix}\right)_{(n-i+1)\times (n-i)},\eqno(2.2)$$
and define
$${\rm rMon}(X) = \left(\begin{smallmatrix}
X_1\oplus\ik\phi_{1}\oplus\cdots\oplus\ik\phi_{n-1}\\
X_{2}\oplus\ik\phi_{2}\oplus\cdots\oplus \ik\phi_{n-1}\\ \vdots
\\X_{n-1}\oplus\ik\phi_{n-1}\\X_n
\end{smallmatrix}\right)_{(\theta_i)}.\eqno(2.3)$$
By construction all $\theta_i$'s are monomorphisms, and hence ${\rm
rMon}(X)\in \mathcal S_n(A)$.

\vskip10pt

\begin{rem} \  By definition ${\rm rMon}(X)$ seems to be
of dependent on the choices of $e_1, \cdots, e_{n-1}$. However, for
an arbitrary choice of $e_1, \cdots, e_{n-1}$, ${\rm rMon}(X)$ will
be proved to be a right minimal approximation of $X$ in $\mathcal
S_n(A)$. Thus, by the uniqueness of a right minimal approximation,
${\rm rMon}(X)$ is in fact independent of the choices of $e_1,
\cdots, e_{n-1}$, up to isomorphism in ${\rm Mor}_n(A)$.
\end{rem}

\subsection{} \ Denote by $\widehat{\mathcal X}$ the full subcategory of
$A$-mod given by ([AR])
$$\widehat{\mathcal X} = \{X \in A\mbox{-}{\rm mod} \
| \ \exists \ \mbox{an exact sequence} \ 0\rightarrow X_{m}
\rightarrow \cdots \rightarrow X_{0}\rightarrow X\rightarrow 0 \
\mbox{with} \  X_{i}\in \mathcal{X}, \ 0\leq i\leq m \}.$$ A
morphism $f: X\rightarrow M$ is {\it right minimal}, if every
endomorphism $g$ of $X$ with $fg = f$ is an isomorphism. {\it A
right approximation of $M$ in $\mathcal X$} is a morphism $f:
X\rightarrow M$ with $X\in\mathcal X$, such that the induced
homomorphism ${\rm{Hom}}_A (X', X)\rightarrow {\rm{Hom}}_A(X', M)$
is surjective for each $X'\in\mathcal X$. A right approximation $f:
X\rightarrow M$ is {\it a right minimal approximation} if $f$ is
right minimal. If every objet $M$ admits a right minimal
approximation in $\mathcal X$, then $\mathcal X$ is called {\it a
contravariantly finite subcategory in} $A$-mod. Dually we have {\it
a covariantly finite subcategory in} $A$-mod. If $\mathcal X$ is
both contravariantly and covariantly finite in $A$-mod, then
$\mathcal X$ is {\it a functorially finite subcategory in} $A$-mod.

\vskip 10pt

\begin{lem} \ Let $A$ be an Artin algebra. Then $\mathcal S_n(A)$ is a
contravariantly finite subcategory in ${\rm Mor}_n(A)$ with
$\widehat{\mathcal S_n(A)} = {\rm Mor}_n(A)$.

Explicitly, for each object $X_{(\phi_i)}$ of \ ${\rm Mor}_n(A)$,
the epimorphism
$$\left(\begin{smallmatrix}(1, 0, \cdots, 0) \\
\vdots\\  (1, 0) \\ 1
\end{smallmatrix}\right): \ \ {\rm rMon}(X) \twoheadrightarrow  X\eqno(2.4)$$
is a right minimal  approximation of $X$ in $\mathcal S_n(A)$.
\end{lem}

\noindent {\bf Proof.} \  By $(2.2)$ and $(2.3)$ it is easy to see
that $(2.4)$ is an epimorphism of ${\rm Mor}_n(A)$. Since $\mathcal
S_n(A)$ is closed under subobjects, it follows that
$\widehat{\mathcal S_n(A)} = {\rm Mor}_n(A)$. Let
$\left(\begin{smallmatrix}
f_1\\
\vdots\\ f_n
\end{smallmatrix}\right): Y_{(\psi_i)} \rightarrow X$ be a morphism of ${\rm Mor}_n(A)$
with $Y = Y_{(\psi_i)}\in \mathcal S_n(A)$. We need to find a
morphism $g = \left(\begin{smallmatrix}
g_1\\
\vdots\\ g_n
\end{smallmatrix}\right): Y \rightarrow {\rm rMon}(X)$
 such that
$$\left(\begin{smallmatrix}
(1, 0, \cdots, 0)\\
\vdots \\ (1, 0) \\ 1
\end{smallmatrix}\right)\left(\begin{smallmatrix}
g_1\\
\vdots\\ g_{n-1}\\ g_n\end{smallmatrix}\right) =
\left(\begin{smallmatrix}
f_1\\
\vdots\\ f_{n-1}\\f_n\end{smallmatrix}\right).\eqno(2.5)$$ We will
inductively
construct $g_i = \left(\begin{smallmatrix} f_i\\ \alpha_{i  i} \\
\vdots \\ \alpha_{i  n-1}
\end{smallmatrix}\right): Y_i \rightarrow X_i\oplus
\ik\phi_{i}\oplus \cdots\oplus\ik\phi_{n-1}, \ 1\le i\le n-1,$ such
that $g: Y \rightarrow {\rm rMon}(X)$ is a morphism of ${\rm
Mor}_n(A)$, i.e., $\theta_i g_{i+1} = g_i\psi_i,$ or explicitly,
such that
$$\left(\begin{smallmatrix} \phi_if_{i+1}\\ e_if_{i+1}
\\ \alpha_{i+1  i+1}\\ \vdots \\ \alpha_{i+1  n-1}
\end{smallmatrix}\right)=\left(\begin{smallmatrix} f_{i}\psi_i\\ \alpha_{i  i}
\psi_i\\ \alpha_{i  i+1} \psi_i\\ \vdots \\ \alpha_{i n-1}\psi_i
\end{smallmatrix}\right).\eqno(2.6)$$

Clearly  $g_n = f_n$. Since $\psi_{n-1}: Y_n\hookrightarrow Y_{n-1}$
is monic and $\ik\phi_{n-1}$ is an injective object, it follows that
the composition $Y_n\stackrel {f_n} \longrightarrow X_n\stackrel
{e_{n-1}} \longrightarrow \ik\phi_{n-1}$ extends to a
morphism $\alpha_{n-1 n-1}: Y_{n-1}\rightarrow \ik\phi_{n-1}$. Define $g_{n-1} = \left(\begin{smallmatrix} f_{n-1}\\
\alpha_{n-1 \ n-1}
\end{smallmatrix}\right)$. Then we have
 \ $\left(\begin{smallmatrix} \phi_{n-1}f_{n}\\ e_{n-1}f_n
\end{smallmatrix}\right)=\left(\begin{smallmatrix} f_{n-1}\psi_{n-1}\\ \alpha_{n-1
n-1} \psi_{n-1}
\end{smallmatrix}\right).$
Assume that we have constructed $g_{n-1}, \cdots, g_t \ (t\ge 2)$,
such that $(2.6)$ holds for $t\le i\le n-1$. Since $\psi_{t-1}:
Y_t\hookrightarrow Y_{t-1}$ is monic and $\ik\phi_{t-1}$ is an
injective object, it follows that the composition $Y_t\stackrel
{f_t} \longrightarrow X_t\stackrel {e_{t-1}} \longrightarrow
\ik\phi_{t-1}$ extends to a morphism $\alpha_{t-1 t-1}:
Y_{t-1}\rightarrow \ik\phi_{t-1}$. Similarly, for $t\le j \le n-1$,
$\alpha_{tj}: Y_t \rightarrow \ik\phi_j$ extends to a morphism
$\alpha_{t-1 j}: Y_{t-1}\rightarrow \ik\phi_j$.
Define $$g_{t-1} = \left(\begin{smallmatrix} f_{t-1}\\
\alpha_{t-1 t-1} \\ \alpha_{t-1 t}\\ \vdots \\ \alpha_{t-1 n-1}
\end{smallmatrix}\right).$$ By construction
$(2.6)$ holds for $i=t-1$, and then $(2.5)$ is clearly satisfied.
This proves that $(2.4)$ is a right approximation of $X$ in
$\mathcal S_n(A)$.

Now we prove that $(2.4)$ is right minimal. Assume that
$\left(\begin{smallmatrix} h_1\\ \vdots\\ h_n
\end{smallmatrix}\right)$ is an endomorphism of ${\rm rMon}(X)$ such
that $\left(\begin{smallmatrix} (1, 0, \cdots, 0)\\ \vdots \\
(1, 0)\\ 1
\end{smallmatrix}\right)\left(\begin{smallmatrix} h_1\\ \vdots\\
h_{n-1}\\ h_n \end{smallmatrix}\right) = \left(\begin{smallmatrix}
(1, 0, \cdots, 0)\\ \vdots \\ (1, 0)\\ 1 \end{smallmatrix}\right).$
We need to prove all $h_i$'s are isomorphisms. Write
$$h_{n-i+1}: \ X_{n-i+1}\oplus \ik\phi_{n-i+1}\oplus.
\cdots\oplus\ik\phi_{n-1}\longrightarrow X_{n-i+1}\oplus
\ik\phi_{n-i+1}\oplus \cdots\oplus\ik\phi_{n-1}$$ as
$\left(\begin{smallmatrix}
h_{n-i+1}^{11} & \cdots &  h_{n-i+1}^{1i} \\
\vdots & \cdots &  \vdots \\ \\ h_{n-i+1}^{i1} & \cdots &
h_{n-i+1}^{ii}
\end{smallmatrix}\right).$
Then $h_{n-i+1}$ is of the form
 \ $h_{n-i+1} = \left(\begin{smallmatrix}
1 & 0 & \cdots &  0 \\ & & & & \\ h_{n-i+1}^{21} & h_{n-i+1}^{22} &\cdots &  h_{n-i+1}^{2i} \\
\vdots & \vdots & \cdots &  \vdots \\ \\ h_{n-i+1}^{i1} &
h_{n-i+1}^{i2} &\cdots & h_{n-i+1}^{ii}
\end{smallmatrix}\right).$
It suffices to prove that all $h_{n-i+1}$ are lower triangular
matrices with diagonal elements being isomorphisms. We do this by
induction. Clearly $h_n = 1$. From the commutative diagram
\[\xymatrix { X_n \ar[r]^-{\theta_{n-1}}\ar@{=}[d] & X_{n-1}\oplus \ik\phi_{n-1} \ar[d]^{h_{n-1}}\\
X_n \ar[r]^-{\theta_{n-1}} & X_{n-1}\oplus \ik\phi_{n-1}}\] we have
 \ $\left(\begin{smallmatrix}
1 & 0 \\ \\ h_{n-1}^{21} & h_{n-1}^{22}
\end{smallmatrix}\right)\left(\begin{smallmatrix}
\phi_{n-1}\\
e_{n-1}
\end{smallmatrix}\right) = \left(\begin{smallmatrix}
\phi_{n-1}\\
e_{n-1}
\end{smallmatrix}\right),$ \
i.e.,  $h_{n-1}^{21}\phi_{n-1} + h_{n-1}^{22}e_{n-1} = e_{n-1}:
X_n\rightarrow \ik\phi_{n-1}$. Restricting the both sides to
$\Ker\phi_{n-1}$ we get $h_{n-1}^{22}e'_{n-1} = e_{n-1}':
\Ker\phi_{n-1} \rightarrow \ik\phi_{n-1}$ (see $(2.1)$). Since
$e_{n-1}'$ is an injective envelope, by definition $h_{n-1}^{22}$ is
an isomorphism.

Assume that $h_{n-t+1} \ (t\ge 2)$ is a lower triangular matrix with
diagonal elements being isomorphisms. Since $\left(\begin{smallmatrix}h_1 \\
\vdots\\  \\ h_n
\end{smallmatrix}\right): \ \ {\rm rMon}(X) \rightarrow {\rm rMon}(X)$ is a morphism,  we have
$$\left(\begin{smallmatrix}
\phi_{n-t}& 0 & 0  &\cdots &  0 \\
e_{n-t}& 0 & 0 & \cdots &  0 \\
0 & 1& 0 & \cdots &  0 \\
0 & 0& 1 & \cdots &  0\\
 \vdots & \vdots & \vdots &
\cdots & \vdots \\
0 & 0& 0 & \cdots & 1
\end{smallmatrix}\right)_{(t+1)\times t}\left(\begin{smallmatrix}
1 & 0 & \cdots &  0 \\ & & & & \\ h_{n-t+1}^{21} & h_{n-t+1}^{22} &\cdots &  h_{n-t+1}^{2t} \\
\vdots & \vdots & \cdots &  \vdots \\ \\ h_{n-t+1}^{t1} &
h_{n-t+1}^{t2} &\cdots & h_{n-t+1}^{tt}
\end{smallmatrix}\right) = \left(\begin{smallmatrix}
1 & 0 & \cdots &  0 \\ & & & & \\ h_{n-t}^{21} & h_{n-t}^{22} &\cdots &  h_{n-t}^{2 t+1} \\
\vdots & \vdots & \cdots &  \vdots \\ \\ h_{n-t}^{{t+1} 1} &
h_{n-t}^{{t+1} 2} &\cdots & h_{n-t}^{{t+1} {t+1}}
\end{smallmatrix}\right)\left(\begin{smallmatrix}
\phi_{n-t}& 0 & 0  &\cdots &  0 \\
e_{n-t}& 0 & 0 & \cdots &  0 \\
0 & 1& 0 & \cdots &  0 \\
0 & 0& 1 & \cdots &  0\\
 \vdots & \vdots & \vdots &
\cdots & \vdots \\
0 & 0& 0 & \cdots & 1
\end{smallmatrix}\right),$$
i.e.,
$$\left(\begin{smallmatrix}
\phi_{n-t}& 0 & 0  &\cdots &  0 \\
e_{n-t}& 0 & 0 & \cdots &  0 \\ \\
h_{n-t+1}^{21} & h_{n-t+1}^{22} & h_{n-t+1}^{23} & \cdots &  h_{n-t+1}^{2t} \\
\vdots & \vdots & \vdots & \cdots & \vdots \\ \\
h_{n-t+1}^{t1} & h_{n-t+1}^{t2} & h_{n-t+1}^{t3} & \cdots &
h_{n-t+1}^{tt}
\end{smallmatrix}\right)_{(t+1)\times t} = \left(\begin{smallmatrix}
\phi_{n-t} & 0 & 0 & \cdots &  0 \\ \\
h_{n-t}^{21}\phi_{n-t}+h^{22}_{n-t}e_{n-t} & h_{n-t}^{23} &
h_{n-t}^{24} &
\cdots &  h_{n-t}^{2 t+1} \\ \\
h_{n-t}^{31}\phi_{n-t}+h^{32}_{n-t}e_{n-t} & h_{n-t}^{33} &
h_{n-t}^{34} &
\cdots &  h_{n-t}^{3 t+1} \\
\vdots & \vdots & \vdots & \cdots &  \vdots \\
\\ h_{n-t}^{t+1 1}\phi_{n-t}+h^{t+1 2}_{n-t}e_{n-t} & h_{n-t}^{t+1 3} & h_{n-t}^{t+1 4} &
\cdots &  h_{n-t}^{t+1 t+1}
\end{smallmatrix}\right)_{(t+1)\times t}.$$
Comparing the second row of the both sides we see
$$h_{n-t}^{23} =0, \ h_{n-t}^{24} = 0, \cdots,
h_{n-t}^{2 t+1} = 0.\eqno(2.7)$$ For $3\le i\le t+1, \ 2\le j\le t$,
comparing the $(i, j)$-entries in both sides we have
$$h_{n-t}^{i j+1} =h_{n-t+1}^{i-1 j}.\eqno(2.8)$$
Since $h_{n-t+1}^{ij} = 0$ for $j>i$ and $h_{n-t+1}^{22}, \cdots,
h_{n-t+1}^{tt}$ are isomorphisms, it follows from $(2.7)$ and
$(2.8)$ that $h_{n-t}^{ij} = h_{n-t+1}^{i-1 j-1} = 0,  \ \forall \
j>i$, and that $h_{n-t}^{ss} = h_{n-t+1}^{s-1 s-1}$ for $s = 3,
\cdots, t+1.$ It remains to prove that $h_{n-t}^{22}$ is an
isomorphism. Comparing the $(2, 1)$-entries we have
$$h_{n-t}^{21}\phi_{n-t}+h^{22}_{n-t}e_{n-t} = e_{n-t}: \  X_{n-t+1}\longrightarrow \ik\phi_{n-t}.$$
Again restricting the both sides to $\Ker\phi_{n-t}$ and by a same
argument we see that $h_{n-t}^{22}$ is an isomorphism. This
completes the proof. $\s$

\vskip10pt

\subsection{\bf Proof of Theorem 2.1.} \ By Corollary 1.4 and Lemma 2.3 $\mathcal S_n(A)$
is a resolving contravariantly finite subcategory in $T_n(A)$-mod.
Then by Corollary 0.3 of [KS] (which asserts that a resolving
contravariantly finite subcategory in $A$-mod is also covariantly
finite in $A$-mod) $\mathcal S_n(A)$ is a functorially finite
subcategory in $T_n(A)$-mod. Thus  $\mathcal S_n(A)$ has
Auslander-Reiten sequences, by Theorem 2.4 of [AS]. $\s$

\vskip10pt

\subsection{} \ For a later use we write down the dual of Theorem 2.1.

\vskip10pt

\noindent {\bf Theorem 2.1'.} \ {\it Let $A$ be an Artin algebra.
Then $\mathcal F_n(A)$ is a functorially finite subcategory in ${\rm
Mor}_n(A)$ and has Auslander-Reiten sequences.}

\vskip10pt

\section{\bf Monomorphism categories and cotilting theory}

The promised reciprocity will be proved, and some consequences will
be given.

\subsection{} \ Let $D$ be the duality $A^{op}\mbox{-}{\rm mod} \rightarrow
A\mbox{-}{\rm mod}$. For $M\in A\mbox{-}{\rm mod}$, denote by ${\rm
add}(M)$ the full subcategory of $A$-mod consisting of all the
direct summands of finite direct sums of copies of $M$, and by $\
^\perp M$ the full subcategory of $A$-mod given by $\{ X\in
A\mbox{-}{\rm mod} \ | \ \Ext^i_A(X, M)= 0, \ \forall \ i\ge 1\}$.

\vskip10pt

An $A$-module $T$ is {\it an $r$-cotilting module} if the following
three conditions are satisfied

\vskip5pt

$(i)$ \ ${\rm inj.dim} T \le r$;

$(ii)$ \ $\Ext^{i}_{A}(T,T)=0$ for $i\geq 1$; and

$(iii)$ \ there is an exact sequence $0\rightarrow T_s \rightarrow
\cdots \rightarrow T_0 \rightarrow D(A_A)\rightarrow 0$ with $T_i\in
{\rm add}(T), \ 0\le i\le s$.

\vskip5pt

\noindent We refer to [HR] and [AR] for the tilting theory.

\vskip10pt

Given an $A$-module $M$, using functor ${\rm \bf m}_i: A\mbox{-}{\rm
mod} \rightarrow T_n(A)\mbox{-}{\rm mod}$ we have a $T_n(A)$-module
$${\rm \bf m}(M) =  \bigoplus\limits_{1\le i\le n} {\rm \bf m}_i(M) =
\left(\begin{smallmatrix}
   M \\
   0\\ 0\\ \vdots \\ 0
\end{smallmatrix}\right) \oplus \left(\begin{smallmatrix}
   M \\
   M\\ 0 \\ \vdots \\ 0
\end{smallmatrix}\right)\oplus \cdots \oplus \left(\begin{smallmatrix}
  M  \\ M \\ M \\ \vdots \\ M
\end{smallmatrix}\right).\eqno(3.1)$$

\vskip10pt

The key result of this paper is as follows.

\vskip10pt

\begin{thm} \ Let $A$ be an Artin algebra, and $T$ an $A$-module.

\vskip10pt

$(i)$ \ If there is an exact sequence $0\rightarrow T_s \rightarrow
\cdots \rightarrow T_0 \rightarrow D(A_A)\rightarrow 0$ with $T_i\in
{\rm add}(T), \ 0\le i\le s$, then \ $\mathcal S_n(^\perp T) = \
^\perp {\rm \bf m}(T).$

\vskip10pt

$(ii)$ \  If $T$ is a cotilting $A$-module, then ${\rm \bf m}(T)$ is
a unique cotilting $T_n(A)$-module, up to multiplicities of
indecomposable direct summands, such that \ $\mathcal S_n(^\perp T)
= \ ^\perp {\rm \bf m}(T).$
\end{thm}

\vskip10pt

Taking $T = D(A_A)$ in Theorem 3.1$(ii)$ we have

\vskip10pt

\begin{cor} \ Let $A$ be an Artin algebra. Then
$\mathcal S_n(A) = \ ^\perp {\rm \bf m}(D(A_A))$.
\end{cor}

In fact, ${\rm \bf m}(D(A_A))$ is the unique cotilting
$T_n(A)$-module, up to multiplicities of indecomposable direct
summands,  such that $\mathcal S_n(A) = \ ^\perp {\rm \bf
m}(D(A_A))$; moreover, ${\rm inj.dim}\ {\rm \bf m}(D(A_A)) =1$, and
$\End_{T_n(A)}({\rm \bf m}(D(A_A))) \cong (T_n(A))^{op}$. Note that
the unique existence of a cotilting $T_n(A)$-module $C$ such that
$\mathcal S_n(A) = \ ^\perp C$ is also guaranteed by Theorem 5.5(a)
in [AR]: since by Lemma 2.3 and Corollary 1.4 $\mathcal S_n(A)$ is a
resolving contravariantly finite subcategory in $T_n(A)$-mod with
$\widehat{\mathcal S_n(A)} = T_n(A)$-mod.

\vskip10pt

If ${\rm inj. dim} A_A < \infty$, we can take $T = \ _AA$ in Theorem
3.1$(i)$ to get

\vskip10pt

\begin{cor} \ Let $A$ be an Artin algebra with ${\rm inj.
dim} A_A < \infty$. Then $\mathcal S_n(^\perp A) = \ ^\perp {\rm \bf
m}(A).$
\end{cor}

\subsection{} \ The proof of Theorem 3.1 needs Theorem 2.1
and the six adjoint pairs between $A$-mod and $T_n(A)$-mod, which
were implied by Lemma 1.2 and will be further explored in the
following.

\vskip10pt

\begin{lem} \ Let $A$ be an Artin algebra and  $M$ an arbitrary $A$-module. Then

\vskip5pt

$(i)$ \ For each $X\in T_n(A)$-mod, we have isomorphisms of abelian
groups, which are natural in both positions
$$\Ext^j_{T_n(A)}({\rm \bf m}_i(M), X) \cong \Ext^j_A (M, X_i), \  \  j\ge 0, \ 1\le i\le n, \eqno(3.2)$$
$$\Ext^j_{T_n(A)}(X, {\rm \bf m}_n(M)) \cong \Ext^j_A(X_1,
M), \ \  j\ge 0, \eqno(3.3)$$
$$\Ext^j_{T_n(A)}(X, {\rm \bf p}_i(M)) \cong \Ext^j_A (X_{n-i+1},
M), \  \  j\ge 0, \ 1\le i\le n, \eqno(3.4)$$
$$\Ext^j_{T_n(A)}({\rm \bf p}_n(M), X) \cong \Ext^j_A (M, X_n), \ j\ge 0. \eqno(3.5)$$

\vskip5pt

$(ii)$ \ For each $X=X_{(\phi_i)}\in \mathcal S_n(A)$, we have
isomorphisms of abelian groups, which are natural in both positions
$$\Ext^j_{T_n(A)}(X, {\rm \bf m}_i(M)) \cong \Ext^j_A (\Cok(\phi_{1}\cdots \phi_{i}), M),
\ \   j\ge 0, \ 1\le i\le n-1. \eqno(3.6)$$

\vskip5pt

$(ii)'$ \ For each $X=X_{(\phi_i)}\in \mathcal F_n(A)$, we have
isomorphisms of abelian groups, which are natural in both positions
$$\Ext^j_{T_n(A)}({\rm \bf p}_i(M), X) \cong
\Ext^j_A(M, \Ker(\phi_{n-i}\cdots \phi_{n-1})), \ \  j\ge 0, \ 1\le
i\le n-1. \eqno(3.7)$$

\end{lem} \noindent{\bf Proof.} \ $(i)$ \ We  justify $(3.3)$.
Taking the $1$-st branch of a projective resolution
$$\cdots \longrightarrow P^1_{(\psi_i^1)}\longrightarrow P^0_{(\psi_i^0)}\longrightarrow
X_{(\phi_i)}\longrightarrow 0\eqno(*)$$ of $X=X_{(\phi_i)}$, by
$(1.8)$ we get a projective resolution $\cdots \rightarrow
P^1_1\rightarrow P^0_1\rightarrow X_1\rightarrow 0$ of $X_1$. On the
other hand, by $(1.4)$ we get the following isomorphic complexes
(for saving the space I omit $\Hom$, and same convention below)
\[\xymatrix{0\ar[r] & (X, {\rm \bf m}_n(M)) \ar[r]\ar[d]^-{\wr}& (P^0, {\rm \bf m}_n(M)) \ar[r]\ar[d]^-{\wr}&
(P^1, {\rm \bf m}_n(M))\ar[r]\ar[d]^-{\wr}&  \cdots \\
0\ar[r] &  (X_1, M) \ar[r] & (P^0_1, M) \ar[r] & (P^1_1, M) \ar[r] &
\cdots.}\] This implies $(3.3)$.

\vskip5pt

$(ii)$ \ By Corollary 1.4 $\mathcal S_n(A)$ is a resolving
subcategory of $T_n(A)$-mod, hence by Lemma 1.1$(i)$ we deduce from
$(*)$ that
$$\cdots \longrightarrow \Cok(\psi_1^1\cdots\psi_i^1)\longrightarrow \Cok(\psi_1^0\cdots\psi_i^0)\longrightarrow
 \Cok(\phi_1\cdots\phi_i) \longrightarrow 0$$
is also exact (it is here we need the assumption $X_{(\phi_i)}\in
\mathcal S_n(A)$).  By $(1.8)$ we see that
$\Cok(\psi_1^j\cdots\psi_i^j)$ is again a projective $A$-module for
every $j$ (it suffices to see
 this for indecomposable projective $T_n(A)$-modules, which is of the form ${\rm\bf m}_i(P)$),
it follows that this exact sequence turns out to be a projective
resolution of $\Cok(\phi_1\cdots\phi_i)$. On the other hand,  for
each $1\le i\le n-1$,  by $(1.3)$ we get the following two
isomorphic complexes
\[\xymatrix{0\ar[r] &(X, {\rm \bf m}_i(M))
\ar[r]\ar[d]^-{\wr}& (P^0, {\rm \bf m}_i(M)) \ar[r]\ar[d]^-{\wr}&
(P^1, {\rm \bf m}_i(M))\ar[r]\ar[d]^-{\wr}&  \cdots \\
0\ar[r] &  (\Cok(\phi_1\cdots\phi_i), M) \ar[r] &
(\Cok(\psi_1^0\cdots\psi_i^0), M) \ar[r] &
(\Cok(\psi_1^1\cdots\psi_i^1), M) \ar[r] & \cdots.}\] This implies
$(3.6)$.   $\s$

\vskip10pt

\subsection{} \ The proof of the following lemma needs Ringel - Schmidmeier's theorem.

\vskip10pt

\begin{lem} \ Let $A$ be an Artin algebra and $X_{(\phi_i)}$ a $T_n(A)$-module.  If $X_{(\phi_i)}\in \ ^\perp {\rm \bf m}(D(A_A))$, then $\phi_i$
is monic, \ $1\le j\le n-1$.
\end{lem}
\noindent {\bf Proof.} \ Taking a right minimal approximation of
$X_{(\phi_i)}$ in $\mathcal S_n(A)$, by $(2.4)$ we have an exact
sequence
$$0\longrightarrow K_{(\theta'_i)}\longrightarrow ({\rm rMon}(X))_{(\theta_i)}
\longrightarrow X_{(\phi_i)}\longrightarrow 0.\eqno(**)$$ Applying
$\Hom_{T_n(A)}(-, {\rm \bf m}_i(D(A_A)))$ to $(**)$ we get an exact
sequence, and by $(1.3)$ this exact sequence is
$$0\longrightarrow (\Cok(\phi_1\cdots\phi_i), D(A_A))\longrightarrow
(\Cok(\theta_1\cdots\theta_i), D(A_A))\longrightarrow
(\Cok(\theta'_1\cdots\theta'_i), D(A_A))\longrightarrow 0,$$ and
hence we get the following exact sequence, which is induced by
$(**)$
$$0\longrightarrow
\Cok(\theta'_1\cdots\theta'_i)\longrightarrow
\Cok(\theta_1\cdots\theta_i)\longrightarrow
\Cok(\phi_1\cdots\phi_i)\longrightarrow 0.$$ On the other hand,
since $\theta_1\cdots\theta_i$ is monic, by Lemma 1.1$(i)$ we get
the following exact sequence, which is again induced by $(**)$
$$0\longrightarrow\Ker(\phi_1\cdots\phi_i)\longrightarrow
\Cok(\theta'_1\cdots\theta'_i)\longrightarrow
\Cok(\theta_1\cdots\theta_i)\longrightarrow
\Cok(\phi_1\cdots\phi_i)\longrightarrow 0.$$ Thus
$\Ker(\phi_1\cdots\phi_i) = 0,$ and hence $\phi_i$ is monic for
$1\le i\le n-1$. $\s$

\vskip10pt Given an $A$-module $M$, using functor ${\rm \bf p}_i:
A\mbox{-}{\rm mod} \rightarrow T_n(A)\mbox{-}{\rm mod}$ we get a
$T_n(A)$-module
$${\rm \bf p}(M) =  \bigoplus\limits_{1\le i\le n} {\rm \bf p}_i(M) =
\left(\begin{smallmatrix}
   0 \\ \vdots \\ 0 \\ 0 \\ M
\end{smallmatrix}\right) \oplus \left(\begin{smallmatrix}
   0 \\
\vdots \\ 0 \\ M\\ M
\end{smallmatrix}\right)\oplus \cdots \oplus \left(\begin{smallmatrix}
  M  \\ M \\ M \\ \vdots \\ M
\end{smallmatrix}\right).$$

\vskip10pt

\begin{prop} \ Let $A$ be an Artin algebra and $T$ an arbitrary
$A$-module.  Then $$\mathcal S_n(^\perp T) = \ ^\perp {\rm \bf m}(T)
\cap \ ^\perp {\rm \bf p}(T)\cap \ ^\perp {\rm \bf m}(D(A_A)).$$
\end{prop}
\noindent {\bf Proof.} \ By $(3.4)$ we have $\mathcal S_n(^\perp T)
\subseteq   \ ^\perp {\rm \bf p}(T).$ By $(3.3)$ and $(3.6)$ we have
$\mathcal S_n(^\perp T) \subseteq   \ ^\perp {\rm \bf m}(D(A_A)).$
Let $X_{(\phi_i)}\in \mathcal S_n(^\perp T)$. By definition $\phi_i:
X_{i+1} \hookrightarrow X_i$ is monic and $\Cok\phi_i\in \ ^\perp T,
\ 1\le i\le n-1.$ By the exact sequence $0\rightarrow \Cok \phi_{i}
\rightarrow \Cok (\phi_{1}\cdots \phi_{i})\rightarrow \Cok
(\phi_{1}\cdots \phi_{i-1})\rightarrow 0$ we inductively see $\Cok
(\phi_{1}\cdots \phi_{i})\in \ ^\perp T$ for $1\le i\le n-1$, and
then by $(3.6)$ and $(3.3)$ this means $X_{(\phi_i)}\in \ ^\perp
{\rm \bf m}(T).$ This proves $\mathcal S_n(^\perp T) \subseteq \
^\perp {\rm \bf m}(T) \cap \ ^\perp {\rm \bf p}(T)\cap \ ^\perp {\rm
\bf m}(D(A_A)).$

Conversely, let $X_{(\phi_i)}\in \ ^\perp {\rm \bf m}(T) \cap \
^\perp {\rm \bf p}(T)\cap \ ^\perp {\rm \bf m}(D(A_A)).$ Then by
$(3.4)$ we have $X_i\in \ ^\perp T, \ 1\le i\le n$, and by Lemma 3.5
$\phi_i: X_{i+1} \hookrightarrow X_i$ is monic, $1\le i\le n-1$. By
$(3.6)$ we know $\Cok (\phi_{1}\cdots \phi_{i})\in \ ^\perp T$ for
$1\le i\le n-1$, and from the exact sequence $0\rightarrow \Cok
\phi_{i} \rightarrow \Cok (\phi_{1}\cdots \phi_{i})\rightarrow \Cok
(\phi_{1}\cdots \phi_{i-1})\rightarrow 0$ we know $\Cok\phi_i\in \
^\perp T, \ 1\le i\le n-1.$ This proves $X_{(\phi_i)}\in \mathcal
S_n(^\perp T)$ and completes the proof. $\s$

\vskip10pt

\subsection{} \ Now we deal with cotilting modules.

\vskip10pt
\begin{lem} \ Let $A$ be an Artin algebra
and $T$ an $r$-cotilting $A$-module. Then ${\rm \bf m}(T)$ is an
$(r+1)$-cotilting $T_n(A)$-module with $\End_{T_n(A)}({\rm \bf
m}(T)) \cong (T_n(\End_A(T)))^{op}$.
\end{lem}

\noindent {\bf Proof.} \ By $(1.2)$ we have $$\Hom_{T_n(A)}({\rm \bf
m}_i(T), {\rm \bf m}_j(T)) \cong
\begin{cases} 0, & i >  j, \\ \End_A(T), & i\le j,\end{cases}$$ we
infer that $\End_{T_n(A)}({\rm \bf m}(T)) \cong
(T_n(\End_A(T)))^{op}$.

Assume that ${\rm inj.dim} T =r$ with a minimal injective resolution
$0\rightarrow T \rightarrow I_0\rightarrow \cdots \rightarrow
I_r\rightarrow 0$. Since ${\rm \bf m}_n: A\mbox{-}{\rm mod}
\rightarrow T_n(A)\mbox{-}{\rm mod}$ is an exact functor and ${\rm
\bf m}_n(I_j) = {\rm \bf p}_n(I_j)$ is an injective $T_n(A)$-module
for each $j$, it follows that ${\rm inj.dim} \ {\rm \bf m}_n(T)= r$.
Similarly, ${\rm inj.dim} \ {\rm \bf p}_{n-i}(T) =r$ for $1\le i\le
n-1$, and then by the following exact sequence
$$0\longrightarrow {\rm \bf m}_i(T) = \left(\begin{smallmatrix}
   T \\ \vdots \\
   T\\ 0 \\ \vdots \\ 0
\end{smallmatrix}\right) \longrightarrow {\rm \bf m}_n(T)= \left(\begin{smallmatrix}
   T \\ \vdots \\
   T\\ T \\ \vdots \\ T
\end{smallmatrix}\right) \longrightarrow {\rm \bf p}_{n-i}(T)=\left(\begin{smallmatrix}
   0 \\
   \vdots \\ 0 \\ T \\ \vdots \\ T
\end{smallmatrix}\right) \longrightarrow 0$$
we see ${\rm inj.dim} \ {\rm \bf m}_i(T) \le r+1$. Thus ${\rm
inj.dim} \ {\rm \bf m}(T) \le r+1$.

\vskip10pt

By $(3.2)$ we have  $$\Ext^j_{T_n(A)}({\rm \bf m}_i(T), {\rm \bf
m}(T)) = \Ext^j_A(T, ({\rm \bf m}(T))_i) = \Ext^j_A(T, \underset
{n-i+1}{\underbrace{T\oplus \cdots \oplus T}}) = 0, \ j\ge 1, \ 1\le
i\le n.$$ This proves $\Ext^j_{T_n(A)}({\rm \bf m}(T), {\rm \bf
m}(T)) = 0$ for $j \ge 1$.

\vskip10pt

Since $T$ is a cotilting $A$-module, we have an exact sequence
$$0\longrightarrow T_s \longrightarrow T_{s-1} \longrightarrow \cdots \longrightarrow T_2\stackrel{d_2}
\longrightarrow T_1\stackrel{d_1}\longrightarrow T_0
\stackrel{d_0}\longrightarrow D(A_A)\longrightarrow 0$$ with every
$T_j\in\add(T)$. Clearly we have an exact sequence
$$0\longrightarrow {\rm \bf m}_n(T_s) \longrightarrow \cdots
 \longrightarrow
{\rm \bf m}_n(T_1)\stackrel{{\rm \bf m}_n(d_1)}\longrightarrow {\rm
\bf m}_n(T_0) \stackrel{{\rm \bf m}_n(d_0)}\longrightarrow {\rm \bf
m}_n(D(A_A)) = {\rm \bf p}_n(D(A_A))\longrightarrow 0$$ with every
${\rm \bf m}_n(T_j)\in\add({\rm \bf m}_n(T)) \subseteq \add ({\rm
\bf m}(T))$. For $1\le i\le n-1$, we have the following exact
sequence of $T_n(A)$-modules
$$0\longrightarrow \left(\begin{smallmatrix}
T_0\\
\vdots \\
T_0 \\ \Ker d_0 \\ \vdots \\ \Ker d_0
\end{smallmatrix}\right)_{(\phi^0_j)}\stackrel{\left(\begin{smallmatrix}
1\\
\vdots \\
1 \\ a \\ \vdots \\ a
\end{smallmatrix}\right)} \longrightarrow {\rm \bf m}_n(T_0)=\left(\begin{smallmatrix}
T_0\\
\vdots \\
T_0 \\ T_0 \\ \vdots \\ T_0
\end{smallmatrix}\right)
\stackrel{\left(\begin{smallmatrix}
0\\
\vdots \\
0 \\ d_0 \\ \vdots \\ d_0
\end{smallmatrix}\right)}\longrightarrow {\rm \bf p}_i(D(A_A))=\left(\begin{smallmatrix} 0\\
\vdots\\ 0 \\ D(A) \\ \vdots \\ D(A)
\end{smallmatrix}\right)\longrightarrow 0,$$
where
$$\phi^0_j = \begin{cases} {\rm id}_{T_0}, & 1\le j\le n-i-1, \\  a, & j=n-i, \\ {\rm id}_{\Ker
d_0}, & n-i+1\le j\le n-1,\end{cases}$$ and $a: \Ker
d_0\hookrightarrow T_0$ is the embedding. Consider the following
sequence of $T_n(A)$-modules
$$0\longrightarrow \left(\begin{smallmatrix}
T_1\\
\vdots \\
T_1 \\ \Ker d_1 \\ \vdots \\ \Ker d_1
\end{smallmatrix}\right)_{(\phi^1_j)}\stackrel f \longrightarrow \left(\begin{smallmatrix}
T_1\\
\vdots \\
T_1 \\ T_1 \\ \vdots \\ T_1
\end{smallmatrix}\right)\oplus \left(\begin{smallmatrix}
T_0\\
\vdots \\
T_0 \\ 0 \\ \vdots \\ 0
\end{smallmatrix}\right)
\stackrel {(\left(\begin{smallmatrix}
d_1\\
\vdots \\
d_1 \\ d_1 \\ \vdots \\ d_1
\end{smallmatrix}\right), \left(\begin{smallmatrix}
1\\
\vdots \\
1 \\ 0 \\ \vdots \\ 0
\end{smallmatrix}\right))} \longrightarrow \left(\begin{smallmatrix}
T_0\\
\vdots \\
T_0 \\ \Ker d_0 \\ \vdots \\ \Ker d_0
\end{smallmatrix}\right)_{(\phi^0_j)}\longrightarrow 0,\eqno(3.8)$$

\vskip10pt\noindent where
$$f =  \left(\begin{smallmatrix}
\left(\begin{smallmatrix} 1\\-d_1
\end{smallmatrix}\right) \\
\vdots \\
\left(\begin{smallmatrix} 1\\-d_1
\end{smallmatrix}\right) \\ a \\ \vdots \\ a
\end{smallmatrix}\right): \ \ \left(\begin{smallmatrix}
T_1\\
\vdots \\
T_1 \\ \\ \\ \Ker d_1 \\ \vdots \\ \Ker d_1
\end{smallmatrix}\right) \longrightarrow \left(\begin{smallmatrix}
T_1\oplus T_0\\
\vdots \\
T_1\oplus T_0 \\ \\ \\ T_1\\ \vdots \\ T_1
\end{smallmatrix}\right).$$
It is routine to see that $f$ and all other maps in $(3.8)$ are
$T_n(A)$-maps (i.e., $(1.1)$ is satisfied for each map). Since both
$0\rightarrow T_1 \stackrel{\left(\begin{smallmatrix}1\\ -d_1
\end{smallmatrix}\right)}\longrightarrow T_1\oplus T_0 \stackrel{(d_1, 1)}\longrightarrow T_0 \rightarrow0$
and $0\rightarrow \Ker d_1 \stackrel{a}\hookrightarrow T_1
\stackrel{d_1}\longrightarrow \Ker d_0 \rightarrow 0 $ are exact, it
follows that $(3.8)$ is exact.

\vskip10pt

Repeating this process we get the following exact sequence of
$T_n(A)$-modules \begin{align*}  0\longrightarrow {\rm \bf
m}_{n-i}(T_s)&\longrightarrow {\rm \bf m}_{n}(T_s)\oplus {\rm \bf
m}_{n-i}(T_{s-1})\longrightarrow \cdots \longrightarrow {\rm \bf
m}_{n}(T_2)\oplus {\rm \bf m}_{n-i}(T_{1})\\ & \longrightarrow {\rm
\bf m}_{n}(T_1)\oplus {\rm \bf m}_{n-i}(T_{0})\longrightarrow {\rm
\bf m}_n(T_0)\longrightarrow {\rm \bf p}_i(D(A_A))\longrightarrow
0.\end{align*} By definition ${\rm \bf m}(T)$ is a cotilting
$T_n(A)$-module. $\s$

\vskip10pt

\subsection{} {\bf Proof of Theorem 3.1.} \ $(i)$ \ By Proposition
3.6 we have $\mathcal S_n(^\perp T) \subseteq \ ^\perp {\rm \bf
m}(T).$ Let $ X=X_{(\phi_i)}\in \ ^\perp {\rm \bf m}(T)$. By the
assumption on $T$ we have an exact sequence
$$0\longrightarrow {\rm \bf m}_i(T_s) \longrightarrow \cdots
\longrightarrow {\rm \bf m}_i(T_0) \longrightarrow {\rm \bf
m}_i(D(A_A))\longrightarrow 0, \ \ 1\le i\le n,$$ with ${\rm \bf
m}_i(T_j)\in \add({\rm \bf m}_i(T)), \ 0\le j\le s$. It follows from
$X_{(\phi_j)}\in \ ^\perp {\rm \bf m}_i(T)$ that $X_{(\phi_j)}\in \
^\perp {\rm \bf m}_i(D(A_A))$, and hence by Lemma 3.5 $\phi_i$ is
monic for $1\le i\le n-1$. Now we can use $(3.6)$ to get $\Cok
(\phi_{1}\cdots \phi_{i})\in \ ^\perp T$ for $1\le i\le n-1$, and by
$(3.3)$ we have $X_1\in \ ^\perp T$.  From the exact sequence
$0\rightarrow X_{i+1}\stackrel{\phi_{1}\cdots \phi_{i}}
\hookrightarrow X_1 \twoheadrightarrow \Cok (\phi_{1}\cdots
\phi_{i}) \rightarrow 0$ we know $X_{i+1}\in \ ^\perp T, \ 1\le i\le
n-1.$ From the exact sequence $0\rightarrow \Cok \phi_{i}
\rightarrow\Cok (\phi_{1}\cdots \phi_{i})\rightarrow \Cok
(\phi_{1}\cdots \phi_{i-1})\rightarrow 0$ we know $\Cok\phi_i\in \
^\perp T, \ 1\le i\le n-1.$ This proves $X_{(\phi_i)}\in \mathcal
S_n(^\perp T)$ and hence $\mathcal S_n(^\perp T) = \ ^\perp {\rm \bf
m}(T).$

\vskip10pt

$(ii)$ \  It follows from $(i)$ and Lemma 3.7 that ${\rm \bf m}(T)$
is a cotilting $T_n(A)$-module with $S_n(^\perp T) = \ ^\perp {\rm
\bf m}(T)$. The remaining uniqueness follows from D. Happel's result
on the number of pairwise non-isomorphic direct summands of a
cotilting module ([H1]). This completes the proof. $\s$

\vskip10pt

\subsection{} \ With the similar arguments one can prove

\vskip10pt

\begin{prop} \ Let $A$ be an Artin algebra and $T$ an
arbitrary $A$-module.  Then $$\mathcal S_n(T^\perp) = {\rm \bf
m}(T)^\perp \cap \ ^\perp {\rm \bf m}(D(A_A)).$$ Moreover, if there
is an exact sequence $0\rightarrow T_s \rightarrow \cdots
\rightarrow T_0 \rightarrow D(A_A)\rightarrow 0$ with $T_i\in {\rm
add}(T), \ 0\le i\le s$, then \ $\mathcal S_n(^\perp T\cap T^\perp)
= \ ^\perp {\rm \bf m}(T)\cap {\rm \bf m}(T)^\perp.$\end{prop}

As an application of Theorem 3.1, we get an answer to the main
problem in Introduction.

\vskip10pt

\begin{thm} \ Let $A$ be an Artin algebra and
$\mathcal X$ a full subcategory of $A$-mod. Then $\mathcal
S_n(\mathcal X)$ is a resolving contravariantly finite subcategory
in $T_n(A)$-mod with $\widehat{\mathcal S_n(\mathcal X)} =
T_n(A)$-mod if and only if $\mathcal X$ is a resolving
contravariantly finite subcategory in $A$-mod with
$\widehat{\mathcal X} = A$-mod.\end{thm} \noindent{\bf Proof.} \ If
$\mathcal X$ is a resolving contravariantly finite subcategory in
$A$-mod with $\widehat{\mathcal X} = A$-mod, then by Theorem 5.5(a)
of Auslander-Reiten [AR] there is a cotilting $A$-module $T$ such
that $\mathcal X = \ ^\perp T$. By Theorem 3.1$(ii)$ ${\rm\bf m}(T)$
is a cotilting $T_n(A)$-module such that $\mathcal S_n(\mathcal X) =
\ ^\perp {\rm\bf m}(T).$ Again by Theorem 5.5(a) in [AR] we know
that $\mathcal S_n(\mathcal X)$ is a resolving contravariantly
finite subcategory in $T_n(A)$-mod with $\widehat{\mathcal
S_n(\mathcal X)} = T_n(A)$-mod. Conversely, assume that $\mathcal
S_n(\mathcal X)$ is a resolving contravariantly finite subcategory
in $T_n(A)$-mod with $\widehat{\mathcal S_n(\mathcal X)} =
T_n(A)$-mod. By Corollary 1.4 $\mathcal X$ is a resolving
subcategory of $A$-mod. Since $\mathcal S_n(\mathcal X)$ is
contravariantly finite in $T_n(A)$-mod with $\widehat{\mathcal
S_n(\mathcal X)} = T_n(A)$-mod, by using functor ${\rm\bf m}_1:
A\mbox{-}{\rm mod} \rightarrow \mathcal S_n(A)$, which induces a
functor ${\rm\bf m}_1: \mathcal X \rightarrow \mathcal S_n(\mathcal
X)$,  we infer that $\mathcal X$ is contravariantly finite
subcategory in $A$-mod with $\widehat{\mathcal X} = A$-mod. $\s$

\vskip10pt

\subsection{} \ For a later use we write down the dual versions of Theorem 3.1,
Corollaries 3.2 and 3.3, Propositions 3.6 and 3.8.

\vskip10pt

\noindent{\bf Theorem 3.1'.} \ {\it Let $A$ be an Artin algebra and
$T$ an arbitrary $A$-module.

\vskip5pt

$(i)$ \ We have $\mathcal F_n(T^\perp) = {\rm \bf p}(T)^\perp\cap
{\rm \bf m}(T)^\perp  \cap  {\rm \bf p}(A)^\perp.$

\vskip5pt

$(ii)$ \ If there is an exact sequence $0\rightarrow A \rightarrow
T_0 \rightarrow \cdots \rightarrow T_s \rightarrow 0$ with $T_i\in
{\rm add}(T), \ 0\le i\le s$, then \ $ \mathcal F_n(T^\perp)= {\rm
\bf p}(T)^\perp.$

\vskip5pt

$(iii)$ \  If $T$ is a tilting $A$-module, then ${\rm \bf p}(T)$ is
a unique tilting $T_n(A)$-module, up to multiplicities of
indecomposable direct summands, such that \ $\mathcal F_n(T^\perp) =
\ {\rm \bf p}(T)^\perp.$

\vskip5pt

$(iv)$ \  ${\rm \bf p}(A)$ is the unique tilting $T_n(A)$-module, up
to multiplicities of indecomposable direct summands,  such that
$\mathcal F_n(A) = \ {\rm \bf p}(A)^\perp$. Moreover, ${\rm
proj.dim}\ {\rm \bf p}(A) =1$, and $\End_{T_n(A)}({\rm \bf p}(A))
\cong (T_n(A))^{op}$.

\vskip5pt

$(v)$  \ If ${\rm inj. dim} _AA < \infty$, then $\mathcal
F_n(D(A_A)^\perp) = {\rm \bf p}(D(A_A))^\perp.$

\vskip5pt

$(vi)$  \ We have $\mathcal F_n(^\perp T) =  \ ^\perp{\rm \bf p}(T)
\cap {\rm \bf p}(A)^\perp.$ \vskip5pt

$(vii)$ \ If there is an exact sequence $0\rightarrow A \rightarrow
T_0 \rightarrow \cdots \rightarrow T_s \rightarrow 0$ with $T_i\in
{\rm add}(T), \ 0\le i\le s$, then \ $\mathcal F_n(T^\perp \cap \
^\perp T) = {\rm \bf p}(T)^\perp \cap  \ ^\perp {\rm \bf p}(T).$}

\subsection{} \ We have the following

\vskip10pt

\begin {rem} \ $(i)$ \ The converse of Theorem
3.1$(i)$ is {\bf not} true. For example, let $k$ be a field and $A$
be the path $k$-algebra of the quiver $1 \bullet \longrightarrow
2\bullet.$ Then $T_2(A)$ is the algebra given by the quiver
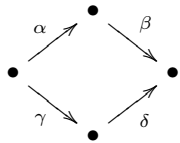
\[\xymatrix @R=1pc  @C=1.5pc{
 & \bullet \ar[rd]^-\beta&\\
 \bullet \ar[ru]^-\alpha\ar[rd]_-\gamma && \bullet \\
 & \bullet\ar[ru]_-\delta
}\]

with relation $\beta\alpha - \delta\gamma$. The Auslander-Reiten
quiver of $T_2(A)$ is

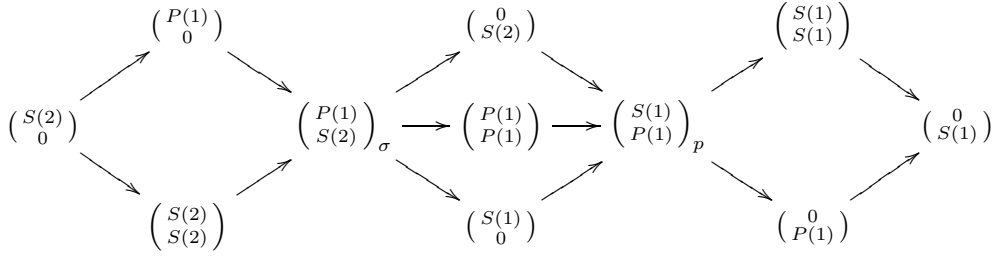
\[\xymatrix @R=1pc  @C=1.5pc{
 &  {\left(\begin{smallmatrix}P(1) \\ 0
\end{smallmatrix}\right)}\ar[rd] & & {\left(\begin{smallmatrix}0 \\ S(2)
\end{smallmatrix}\right)} \ar[rd] & & {\left(\begin{smallmatrix}S(1) \\ S(1)
\end{smallmatrix}\right)} \ar[rd]\\
 {\left(\begin{smallmatrix}S(2) \\ 0
\end{smallmatrix}\right)} \ar[ru]\ar[rd]&& {\left(\begin{smallmatrix}P(1) \\ S(2)
\end{smallmatrix}\right)_\sigma} \ar[ru]\ar[r]\ar[rd] & {\left(\begin{smallmatrix}P(1) \\ P(1)
\end{smallmatrix}\right)} \ar[r]&{\left(\begin{smallmatrix}S(1) \\ P(1)
\end{smallmatrix}\right)_p}
 \ar[ru]\ar[rd] && {\left(\begin{smallmatrix}0 \\ S(1)
\end{smallmatrix}\right)} \\
 & {\left(\begin{smallmatrix}S(2) \\ S(2)
\end{smallmatrix}\right)}\ar[ru] & & {\left(\begin{smallmatrix}S(1) \\ 0
\end{smallmatrix}\right)} \ar[ru] &&{\left(\begin{smallmatrix}0 \\ P(1)
\end{smallmatrix}\right)} \ar[ru]
}\]

\vskip10pt

\noindent Let $T = S(1)\oplus S(2)$. Then $\mathcal S_2(^\perp T) =
\mathcal S_2(^\perp S(2)) = \mathcal S_2({\rm add} (A))$, and
$$^\perp {\rm
\bf m}(T) = \ ^\perp(\left(\begin{smallmatrix}S(1)\oplus S(2) \\ 0
\end{smallmatrix}\right) \oplus \ \left(\begin{smallmatrix}S(1)\oplus S(2) \\ S(1)\oplus S(2)
\end{smallmatrix}\right)) = \ ^\perp\left(\begin{smallmatrix}S(1) \\ 0
\end{smallmatrix}\right)\cap \ ^\perp\left(\begin{smallmatrix}S(2) \\ 0
\end{smallmatrix}\right)  \cap \  ^\perp\left(\begin{smallmatrix}S(2) \\ S(2)
\end{smallmatrix}\right).$$
It is clear that
$$\mathcal S_2(^\perp T) = \mathcal S_2({\rm add} (A)) = {\rm add} (\left(\begin{smallmatrix}S(2) \\ 0
\end{smallmatrix}\right)\oplus \left(\begin{smallmatrix}P(1) \\ 0
\end{smallmatrix}\right)\oplus \left(\begin{smallmatrix}S(2) \\ S(2)
\end{smallmatrix}\right)\oplus \left(\begin{smallmatrix}P(1) \\ P(1)
\end{smallmatrix}\right))
= \  ^\perp {\rm \bf m}(T),$$ but $T$ does not satisfy the condition
in Theorem 3.1$(i)$.

\vskip10pt

Nevertheless, even in this example, for many $A$-modules $T$ {\it
not} satisfying the condition in Theorem 3.1$(i)$,  we have
$\mathcal S_n(^\perp T) \ne \ ^\perp {\rm \bf m}(T)$. For examples,
this is the case when $T = S(1)$, or $T=S(2)$, or $T= P(1)$.

\vskip10pt

$(ii)$  \ Many cotilting $T_n(A)$-modules are {\it not} of the form
${\rm \bf m}(T)$, where $T$ is a cotilting $A$-module. For example,
if $k$ is a field, then $T_3(k)$ is the path $k$-algebra of the
quiver $1 \bullet \longrightarrow 2 \bullet \longrightarrow
3\bullet.$ There two basic cotilting $T_3(k)$-modules having the
simple module $S(2)$ as a direct summand, which are not of the form
${\rm \bf m}(T)$, where $T\in k$-mod.\end{rem}
\section{\bf Application to Gorenstein algebras}

Applying Theorem 3.1 to Gorenstein algebras, we explicitly determine
all the Gorenstein-projective $T_n(A)$-modules. We characterize
self-injective algebras by monomorphism categories.

\vskip10pt

\subsection{} \ Modules in \ $^\perp (_AA)$ are called {\it Cohen-Macaulay $A$-modules}. Denote
$^\perp A$ by ${\rm \bf CM}(A)$. An $A$-module $G$ is {\it
Gorenstein-projective}, if there is an exact sequence $\cdots
\rightarrow P^{-1}\rightarrow P^{0} \stackrel{d^0}{\rightarrow}
P^{1}\rightarrow \cdots$ of projective $A$-modules, which stays
exact under ${\rm Hom}_A(-, A)$, and such that $G\cong
\operatorname{Ker}d^0$. Let $A\mbox{-}\mathcal Gproj$ be the full
subcategory of $A$-mod
 of Gorenstein-projective modules. Then $A\mbox{-}\mathcal
Gproj \subseteq \ {\rm \bf CM}(A)$;  and if $A$ is {\it a Gorenstein
algebra} (i.e., ${\rm inj.dim}\ _AA< \infty$ and ${\rm inj.dim} \
A_A < \infty$), then \ $A\mbox{-}\mathcal Gproj = {\rm \bf CM}(A)$
(Enochs - Jenda [EJ2], Corollary 11.5.3). Determining all the
Cohen-Macaulay $A$-modules and all the Gorenstein-projective
$A$-modules in explicit way, is a basic requirement in applications
(see e.g. [AM], [B], [BGS], [CPST], [EJ2], [GZ], [K]).

\vskip10pt

\begin{cor} \  $(i)$ \ Let $A$ be an
Artin algebra with ${\rm inj.dim} A_A < \infty$. Then $${\rm \bf
CM}(T_n(A)) = \mathcal S_n({\rm \bf CM}(A)).$$

$(ii)$ \  Let $A$ be a Gorenstein algebra. Then
$T_n(A)\mbox{-}\mathcal Gproj = \mathcal S_n(A\mbox{-}\mathcal
Gproj).$
\end{cor}

\noindent {\bf Proof.} \ $(i)$ is a reformulation of Corollary 3.3
since ${\rm \bf m}(_AA) = \ _{T_n(A)}T_n(A)$. If $A$ is Gorenstein,
then it is well-known that $T_n(A)$ is again Gorenstein (for $n=2$
see e.g. [FGR] or [H2]; in general see e.g. [XZ], Lemma 4.1$(i)$),
and hence $(ii)$ follows from $(i)$. $\s$

\vskip10pt Corollary 4.1$(ii)$ was obtained for $n=2$ in Theorem
1.1$(i)$ of [LZ2] (see also Proposition 3.6$(i)$ of [IKM]).

\vskip10pt

\subsection {} \ Dually, denote $D(A_A)^\perp$ by ${\rm
\bf CoCM}(A)$. An $A$-module $G$ is {\it Gorenstein-injective}
([EJ1]), if there is an exact sequence $\cdots \rightarrow
I^{-1}\rightarrow I^{0} \stackrel{d^0}{\rightarrow} I^{1}\rightarrow
\cdots$ of injective $A$-modules, which stays exact under ${\rm
Hom}_A(D(A_A), -)$, and such that $G\cong \operatorname{Ker}d^0$.
Let $A\mbox{-}\mathcal Ginj$ be the full subcategory of $A$-mod of
Gorenstein-injective modules. Then $A\mbox{-}\mathcal Ginj \subseteq
\ {\rm \bf CoCM}(A)$; and if $A$ is Gorenstein then \
$A\mbox{-}\mathcal Ginj = {\rm \bf CoCM}(A)$. By Theorem 3.1'$(v)$
and Corollary 4.1 we have

\vskip10pt

\begin{cor} \  $(i)$ \ Let $A$ be an
Artin algebra with ${\rm inj.dim} _AA < \infty$. Then ${\rm \bf
CoCM}(T_n(A)) = \mathcal F_n({\rm \bf CoCM}(A)).$

\vskip5pt

$(ii)$ \  Let $A$ be a Gorenstein algebra. Then
$T_n(A)\mbox{-}\mathcal Ginj = \mathcal F_n(A\mbox{-}\mathcal
Ginj),$ and the set of $T_n(A)$-modules which are simultaneously
Gorenstein-projective and Gorenstein-injective is
$$\{ {\rm\bf m}_n(M) \ | \ M \ \mbox{is simultaneously a Gorenstein-projective and
Gorenstein-injective A-module}\}.$$
\end{cor}

\vskip10pt

\subsection {} \ Let $D^b(A)$ be the bounded derived category of $A$, and
$K^b(\mathcal P(A))$ the bounded homotopy category of $\mathcal
P(A)$. The singularity category $D^b_{sg}(A)$ of $A$ is defined to
be the Verdier quotient $D^b(A)/K^b(\mathcal P(A))$. If $A$ is
Gorenstein, then there is a triangle-equivalence $D^b_{sg}(A)\cong
\underline {{\rm \bf CM}(A)},$ where $\underline {{\rm \bf CM}(A)}$
is the stable category of \ ${\rm \bf CM}(A)$ modulo $\mathcal P(A)$
([H2], Theorem 4.6; see also [Buc], Theorem 4.4.1). Thus by
Corollary 4.1 we have

\vskip10pt

\begin{cor} \  Let $A$ be a
Gorenstein algebra. Then there is a triangle-equivalence \
$$D^b_{sg}(T_n(A))\cong  \underline {\mathcal S_n({\rm \bf CM}(A))}.$$

In particular, if $A$ is a self-injective algebra, then
$D^b_{sg}(T_n(A))\cong \underline {\mathcal S_n(A)}.$ \end{cor}

\vskip10pt

\subsection {} \ We have the following characterization of self-injective
algebras.

\vskip10pt

\begin{thm} \ Let $A$ be an Artin algebra. Then $A$ is a
self-injective algebra if and only if \ $T_n(A)\mbox{-}\mathcal
Gproj = \mathcal S_n(A).$\end{thm} \noindent {\bf Proof.} \ The
``only if" part follows from Corollary 4.1$(ii)$. Conversely, by
assumption $\left(\begin{smallmatrix}
   D(A_A) \\
   0\\ \vdots \\  0
\end{smallmatrix}\right)\in \mathcal S_n(A)$ is a Gorenstein-projective $T_n(A)$-module. Then
there is an exact sequence $\cdots \rightarrow
P^{-1}\rightarrow P^{0} \stackrel{d^0}{\rightarrow} P^{1}\rightarrow
\cdots $ of projective $T_n(A)$-modules with
$\left(\begin{smallmatrix}
   D(A_A) \\
   0\\ \vdots \\  0
\end{smallmatrix}\right)\cong \operatorname{Ker}d^0$.
By taking the $1$-st branch we get an exact sequence $\cdots
\rightarrow P^{-1}_1\rightarrow P^{0}_1
\stackrel{d^0_1}{\rightarrow} P^{1}_1\rightarrow \cdots $ of
projective $A$-modules with $\Ker d^0_1\cong D(A_A)$. This implies
that $D(A_A)$ is a projective module, i.e., $A$ is self-injective.
$\s$

\section{\bf Finiteness of monomorphism categories}

This section is to characterize $\mathcal S_n(A)$ which is of finite
type.

\vskip10pt

\subsection{} \  An additive full subcategory $\mathcal X$ of $A$-${\rm
mod}$, which is closed under direct summands, is of {\it finite
type} \ if there are only finitely many isomorphism classes of
indecomposable $A$-modules in $\mathcal X$. If $A\mbox{-}\mathcal
Gproj$ is of finite type, then $A$ is said to be {\it CM-finite}.
\vskip10pt

An $A$-module $M$ is  {\it an $A$-generator} if each projective
$A$-module is in $\add M$. A $T_n(A)$-generator $M$ is {\it a
bi-generator of \ $\mathcal S_n(A)$}  if $M\in \mathcal S_n(A)$ and
${\rm \bf m}(D(A_A))\in \add(M)$.

\vskip10pt

\begin{thm} \ Let $A$ be an Artin algebra.
Then $\mathcal S_n(A)$ is of finite type if and only if there is a
bi-generator $M$ of $\mathcal S_n(A)$ such that $\gld
\End_{T_n(A)}(M) \leq 2$.
\end{thm}

\vskip10pt

If $A$ is self-injective, then ${\rm \bf m}(D(A_A))$ is a projective
$T_n(A)$-module, and hence in $\add(M)$ for each $T_n(A)$-generator
$M$.  Combining Theorems 5.1 and 4.4 we have

\vskip10pt

\begin{cor} \ Let $A$ be a self-injective algebra. Then $T_n(A)$ is CM-finite if and only if there
is a $T_n(A)$-generator $M$ which is Gorenstein-projective, such
that $\gld {\rm End}_{T_n(A)}(M) \leq 2$.\end{cor}

\vskip10pt

Corollary 5.2 also simplifies the result in [LZ1] in this special
case. \vskip10pt

\subsection{} \ The proof of Theorem 5.1 will use Corollary 3.2, and Auslander's
idea in proving his classical result cited in Introduction. Given
modules $M, X \in A\mbox{-mod}$, denote by $\Omega_M(X)$ the kernel
of a minimal right approximation $M' \rightarrow X$ of $X$ in $\add
(M)$. Define $\Omega_M^0(X) = X$, and $\Omega_M^{i}(X) =
\Omega_M(\Omega_M^{i-1}(X))$ for $i\ge 1$. Define ${\rm rel.dim}_M
X$ to be the minimal non-negative integer $d$ such that
$\Omega_M^d(X)\in\add (M)$, or $\infty$ if otherwise. The following
fact is well-known.

\vskip10pt

\begin{lem} \ (M. Auslander) \ Let $A$ be an Artin algebra and
$M$ be an $A$-module, and $\Gamma = ({\rm End}_A(M))^{op}.$ Then
${\rm proj.dim}\ _\Gamma\Hom_A(M, X) \le {\rm rel.dim}_M X$ for all
$A$-modules $X$. If $M$ is a generator, then equality holds.
\end{lem}

\vskip10pt

For an $A$-module $T$, denote by $\mathcal{X}_{_{T}}$ the full
subcategory of $A$-mod given by
$$\{X \ | \ \exists \
\mbox{an exact sequence} \ 0\rightarrow X \rightarrow
T_0\stackrel{d_{0}}\rightarrow T_{1}\stackrel{d_{1}}\rightarrow
\cdots, \mbox{with} \ T_{i}\in {\rm add}(T), \ \Ker d_{i}\in
{^{\perp}T, \ \forall \ i\geq 0}\}.$$ Note that
$\mathcal{X}_{T}\subseteq \ ^\perp T$, and $\mathcal{X}_{T}= \
^\perp T$ if $T$ is a cotilting module ([AR], Theorem 5.4$(b)$).

\vskip10pt

\begin {lem} \ Let $A$ be an Artin algebra and
$M$ be an $A$-generator with \ $\Gamma = ({\rm End}_A(M))^{op},$ and
$T\in \add(M)$. Then for every $A$-module $X\in \mathcal{X}_{T}$ and
$X\notin \add(T)$, there is a $\Gamma$-module $Y$ such that ${\rm
proj.dim}_\Gamma Y= 2 + {\rm proj.dim}_\Gamma \Hom_A(M,
X).$\end{lem}

\noindent {\bf Proof.} \ By $X\in \mathcal{X}_{T}$ there is an exact
sequence $0 \rightarrow X \stackrel{u}\rightarrow T_0
\stackrel{v}\rightarrow T_1$ with $T_0,  \ T_1 \in \add(T)\subseteq
\add(M)$. This yields an exact sequence
$$0\longrightarrow \Hom_A(M, X) \stackrel{u_*}\longrightarrow \Hom_A(M,
T_0) \stackrel{v_*}\longrightarrow \Hom_A(M, T_1) \longrightarrow
\Cok v_*\longrightarrow 0.$$ Since the image of $v_*$ is not
projective (otherwise, $u_*$ splits, and then one can deduce a
contradiction $X\in \add(T)$), putting $Y=\Cok v_*$ we have ${\rm
proj.dim}_\Gamma Y = 2+ {\rm proj.dim}_\Gamma \Hom_A(M, X).$ \hfill
$\s$
\subsection {\bf Proof of Theorem 5.1.} \ Assume that $\mathcal S_n(A)$ is of finite type.
Then there is a $T_n(A)$-module $M$ such that $\mathcal S_n(A) =
\add(M)$. Since ${\rm \bf m}(D(A_A))\in \mathcal S_n(A) = \add(M)$
and $\mathcal S_n(A)$ contains all the projective $T_n(A)$-modules,
by definition $M$ is a bi-generator of $\mathcal S_n(A)$. Put
$\Gamma = (\End_{T_n(A)}(M))^{op}$. For every $\Gamma$-module $Y$,
take a projective presentation $
 \Hom_{T_n(A)}(M, M_1) \stackrel{f_*}\longrightarrow \Hom_{T_n(A)}(M, M_0) \rightarrow Y \rightarrow 0
$ of $Y$, where $M_1, M_0 \in \add (M),$ and $f: M_1\rightarrow M_0$
is a $T_n(A)$-map. Since  ${\rm inj.dim} \ {\rm \bf m}(D(A_A))$ $ =
1$ (Lemma 3.7) and $M_1\in \mathcal S_n(A) = \ ^\perp {\rm \bf
m}(D(A_A))$ (Corollary 3.2), it follows that $\Ker f\in \ ^{\perp}
({\rm \bf m}(D(A_A))) = \add(M)$. Thus
$$0\longrightarrow \Hom_{T_n(A)}(M, \Ker f)\longrightarrow \Hom_{T_n(A)}(M, M_1)\longrightarrow \Hom_{T_n(A)}(M, M_0)\longrightarrow
Y\longrightarrow 0$$ is a projective resolution of $\Gamma$-module
$Y$, i.e., ${\rm proj.dim} _\Gamma Y \le 2.$ This proves $\gld
\Gamma \leq 2$. Since $A$ is an Artin algebra, we have $\gld
\End_{T_n(A)}(M) = \gld \Gamma \le 2$.

Conversely, assume that there is a bi-generator $M$ of $\mathcal
S_n(A)$ such that $\gld \End_{T_n(A)}(M) \le 2$. Put $\Gamma =
(\End_{T_n(A)}(M))^{op}$. Then $\gld \Gamma \leq 2$. We claim that
$\add(M) = \ ^\perp {\rm \bf m}(D(A_A))$, and hence by Corollary 3.2
$\mathcal S_n(A)$ is of finite type. In fact, since $M\in \mathcal
S_n(A) = \ ^\perp {\rm \bf m}(D(A_A))$, it follows that $\add(M)
\subseteq \ ^\perp {\rm \bf m}(D(A_A))$. On the other hand, let
$X\in \ ^\perp {\rm \bf m}(D(A_A))$. By Corollary 3.2 ${\rm \bf
m}(D(A_A))$ is a cotilting $T_n(A)$-module, and hence $^\perp {\rm
\bf m}(D(A_A)) = \mathcal X_{{\rm \bf m}(D(A_A))}$, by Theorem
5.4$(b)$ in [AR]. If $X\in\add({\rm \bf m}(D(A_A)))$, then $X\in
\add(M)$ since by assumption ${\rm\bf m}(D(A_A))\in\add(M)$. If
$X\notin\add({\rm \bf m}(D(A_A)))$, then by Lemma 5.4 there is a
$\Gamma$-module $Y$ such that ${\rm proj.dim}_\Gamma Y = 2 + {\rm
proj.dim}_\Gamma \Hom_{T_n(A)}(M, X).$ Now by Lemma 5.3 we have
$${\rm rel.dim}_M X = {\rm proj.dim}_\Gamma\Hom_{T_n(A)}(M, X) =
{\rm proj.dim}_\Gamma Y - 2\le \gld \Gamma -2 \le 0,$$ this means
$X\in\add (M)$. This proves the claim and completes the proof. $\s$


\vskip10pt


\begin{thebibliography}{99}


\bibitem[Ar]{Ar} D. M. Arnold, Abelian groups and representations of finite partially
ordered sets, Canad. Math. Soc. Books in Math., Springer-Verlag, New
York, 2000.
\bibitem[Au]{Au} M. Auslander, Representation dimension of
artin algebras, Queen Mary College Math. Notes, London, 1971. Also
in: Selected works of Maurice Auslander, Part 1 (II),  Edited by I.
Reiten, S. Smal${\o}$, ${\O}$. Solberg, Amer. Math. Soc. (1999),
505-574.
\bibitem[AR]{AR} M. Auslander, I. Reiten, Applications of
contravariantly finite subcategories, Adv. Math. 86(1991), 111-152.
\bibitem[ARS]{ARS} M. Auslander, I. Reiten, S. O.
Smal${\o}$, Representation Theory of Artin Algebras, Cambridge
Studies in Adv. Math. 36., Cambridge Univ. Press,
1995.
\bibitem[AS]{AS} M. Auslander, S. O. Smal${\o}$, Almost split
sequences in subcategories, J. Algebra 69(1981), 426-454.
\bibitem[AM] {AM} L. L. Avramov, A. Martsinkovsky, Absolute, relative, and Tate cohomology of modules of
finite Gorenstein dimension, Proc. London Math. Soc. 85(3)(2002),
393-440.

\bibitem[B]{B} A.
Beligiannis, Cohen-Macaulay modules, (co)tosion pairs and virtually
Gorenstein algebras, J. Algebra 288(1)(2005), 137-211.
\bibitem[Bir]{Bir} G.
Birkhoff, Subgroups of abelian groups, Proc. Lond. Math. Soc. II,
Ser. 38(1934), 385-401.
\bibitem[Buc]{B} R.-O. Buchweitz, Maximal Cohen-Macaulay modules and Tate cohomology over Gorenstein rings,
Unpublished manuscript, Hamburg (1987), 155pp.
\bibitem[BGS]{BGS}
R.-O. Buchweitz, G.-M. Greuel, F.-O. Schreyer, Cohen-Macaulay
modules on hypersurface singularities II, Invent. Math. 88(1)(1987),
165-182.
\bibitem[C]{C}
X. W. Chen, Stable monomorphism category of Frobenius category,
avaible in arXiv: 0911.1987.

\bibitem[CPST]{CPST} L. W. Christensen, G.
Piepmeyer, J. Striuli, R. Takahashi, Finite Gorenstein
representation type implies simple singularity,  Adv. Math.
218(2008), 1012-1026.
\bibitem[EJ1]{EJ1}
E. E. Enochs, O. M. G. Jenda, Gorenstein injective and projective
modules, Math. Z. 220(4)(1995), 611-633.
\bibitem[EJ2]{EJ2} E. E. Enochs, O. M. G. Jenda, Relative homological algebra, De Gruyter Exp. Math. 30.
Walter De Gruyter Co., 2000.
\bibitem[FGR]{FGR} R. Fossum, P.
Griffith, I. Reiten, Trivial extensions of abelian categories,
Lecture Notes in Math. 456, Springer-Verlag, 1975.
\bibitem[GZ]{GZ}
N. Gao, P. Zhang, Gorenstein derived categories, J. Algebra
323(2010), 2041-2057.

\bibitem[H1]{H1} D. Happel, Triangulated categories in
representation theory of finite dimensional algebras, London Math.
Soc. Lecture Notes Ser. 119, Cambridge Uni. Press, 1988.
\bibitem[H2]{H2} D. Happel, On Gorenstein algebras, in:
Representation theory of finite groups and finite-dimensional
algebras, Prog. Math. 95, 389-404, Birkh\"user, Basel, 1991.
\bibitem[HR]{H1} D. Happel, C. M. Ringel, Tilted algebras,
Trans. Amer. Math. Soc. 274(2)(1982), 399-443.
\bibitem[IKM]{IKM} O.
Iyama,  K. Kato, J. I. Miyachi, Recollement on homotopy categories
and Cohen-Macaulay modules, avaible in arXiv: math. RA 0911.0172.

\bibitem[K]{K} H. Kn\"orrer, Cohen-Macaulay
modules on hypersurface singularities I, Invent. Math. 88(1)(1987),
153-164.
\bibitem[KS]{KS} H. Krause, $\O$. Solberg, Applications of cotorsion
pairs,  J. London Math. Soc. 68(3)(2003), 631-650.
\bibitem[KLM]{KLM}
D. Kussin, H. Lenzing, H. Meltzer, Nilpotent operators and weighted
projective lines, avaible in arXiv: math. RT 1002.3797.
\bibitem[LZ1]{LZ1} Z. W. Li, P. Zhang, Gorenstein algebras
of finite Cohen-Macaulay type, Adv. Math. 223(2010), 728-734.
\bibitem[LZ2]{LZ2} Z. W. Li, P. Zhang, A construction of Gorenstein-projective modules, J.
Algebra 323(2010),  1802-1812.
\bibitem[RW]{RW} F. Richman, E. A. Walker, Subgroups of $p^5$-bounded groups, in:
Abelian groups and modules, Trends Math., Birkh\"auser, Basel, 1999,
55-73.
\bibitem[RS1]{RS1} C. M. Ringel, M. Schmidmeier, Submodules
categories of wild representation type, J. Pure Appl. Algebra
205(2)(2006), 412-422.
\bibitem[RS2]{RS2}
C. M. Ringel, M. Schmidmeier, The Auslander-Reiten translation in
submodule categories, Trans. Amer. Math. Soc. 360(2)(2008), 691-716.
\bibitem[RS3]{RS3}
C. M. Ringel, M. Schmidmeier, Invariant subspaces of nilpotent
operators I, J. rein angew. Math. 614 (2008), 1-52.
\bibitem[S]{S} D. Simson, Representation types of the category of subprojective
representations of a finite poset over $K[t]/(t^m)$ and a solution
of a Birkhoff type problem, J. Algebra 311(2007), 1-30.
\bibitem [SW]{SW} D. Simson,  M. Wojewodzki, An algorithmic solution of a
Birkhoff type problem, Fundamenta Informaticae 83(2008), 389-410.
\bibitem[XZ]{XZ} B. L. Xiong,
P. Zhang, Cohen-Macaulay modules over triangular matrix Artin
algebras, preprint (2009).

\end{thebibliography}
\end{document}